\documentclass[letterpaper, 10pt, conference]{ieeeconf}

\IEEEoverridecommandlockouts
\usepackage[T1]{fontenc}
\usepackage[utf8]{inputenc}
\usepackage{amsmath,amssymb,amsfonts,bm}
\usepackage{graphicx}
\usepackage{xcolor}
\usepackage{booktabs}
\usepackage{array}
\usepackage{cite}
\usepackage{tikz}
\usepackage[linesnumbered,ruled,vlined]{algorithm2e}
\SetKwInput{KwIn}{Inputs}
\usepackage{mathtools}
\usepackage{amsthm}
\usepackage[caption=false,font=footnotesize]{subfig}
\usepackage[hidelinks]{hyperref}
\allowdisplaybreaks

\usetikzlibrary{arrows.meta,backgrounds}
\usepackage{balance}

\makeatletter
\renewcommand{\section}{%
  \@startsection{section}{1}{\z@}%
  {1.6ex plus 0.6ex minus 0.3ex}%
  {0.3ex plus 0.2ex}%
  {\normalfont\normalsize\centering\scshape}}

\renewcommand{\subsection}{%
  \@startsection{subsection}{2}{\z@}%
  {1.3ex plus 0.5ex minus 0.3ex}%
  {0.2ex plus 0.2ex}%
  {\normalfont\normalsize\itshape}}
\makeatother

\AtBeginDocument{
  \setlength{\abovedisplayskip}{4pt plus 1pt minus 1pt}
  \setlength{\belowdisplayskip}{4pt plus 1pt minus 1pt}
  \setlength{\abovedisplayshortskip}{2pt plus 1pt}
  \setlength{\belowdisplayshortskip}{3pt plus 1pt minus 1pt}
  \setlength{\jot}{2pt}
}

\newtheorem{assumption}{Assumption}
\newtheorem{theorem}{Theorem}

\theoremstyle{plain}
\newtheorem{lemma}[theorem]{Lemma}

\definecolor{lightgreen}{RGB}{225, 245, 230}
\title{\textbf{Stochastic Optimization with Network Distribution Dynamics \\ via Local Jacobian Approximation}}

\author{
Yuchen Cai, Zhiyu He, Michael Muehlebach, and Florian Dörfler%
\thanks{Yuchen Cai, Zhiyu He, and Florian Dörfler are with the Automatic Control Laboratory, ETH Z\"urich, Switzerland (email: \{yuccai,~zhiyhe,~dorfler\}@ethz.ch). Michael Muehlebach is with the WU Vienna University of Economics and Business, Austria (email: michael.muehlebach@wu.ac.at). This work was supported by the Max Planck ETH Center for Learning Systems and SNSF via NCCR Automation (51NF40\_225155).}
}
\begin{document}

\maketitle

\begin{abstract}
Network systems feature complex interactions among agents, whereby local decisions can influence the state distributions observed by other agents through coupled network dynamics.  Such endogenous distribution shifts across large-scale networks pose significant challenges and call for distributed and adaptive decision-making pipelines.
To address these challenges and collectively optimize a steady-state performance objective, we propose a distributed network-aware decision-response (D-NAD) gradient algorithm. This algorithm enables each agent to approximate network-response sensitivity using only local information (i.e., the Jacobians of local dynamics and objectives) and neighborhood communication.
% These distribution shifts affect the steady-state objective gradient through the network-response sensitivity, which is determined by the Jacobians of the coupled dynamics.
We establish a non-asymptotic closed-loop convergence guarantee and quantify the impact of the Jacobian approximation errors therein. Numerical experiments on a semi-synthetic networked electric-vehicle price--demand problem demonstrate the effectiveness of the proposed algorithm.

\end{abstract}

% \begin{IEEEkeywords}
% Decision-dependent stochastic optimization, network systems, distributed algorithm.
% \end{IEEEkeywords}

\section{Introduction}
% \zhmargin{to be revised}
% Stochastic optimization typically assumes that uncertainty follows a
% fixed distribution independent of the decision. This assumption fails
% when deployed decisions change future observations, as when prices
% affect both revenue and product demand. In such settings, a decision
% changes not only the current objective but also the distribution under
% which future decisions are evaluated. Such feedback is studied in
% performative prediction and decision-dependent stochastic optimization
% \cite{perdomo2020performative,mendler2020stochastic,
% drusvyatskiy2023stochastic,he2025decisiondependent}.

There is growing interest in the feedback loop between the deployed decisions and the underlying data-generating distribution~\cite{perdomo2020performative,hardt2025performative}. Such a feedback phenomenon arises in many practical scenarios, including  strategic classification~\cite{mendler2020stochastic,drusvyatskiy2023stochastic}, recommendation systems~\cite{he2025decisiondependent,wu2026reaching}, and dynamic pricing~\cite{cavraro2026demandresponse,wang2025constrainedddd}, where players with market power impact the demand distribution through their pricing decisions. In these settings, decision-makers (a centralized entity or multiple interdependent agents) cause endogenous distribution shifts, which in turn affect future problems and subsequent decisions. 
% Many of these applications have a multi-agent structure, as illustrated by our electric-vehicle (EV) price--demand case study in Section~\ref{sec:EV_model}, where a local price change can affect adoption across coupled markets.

The above decision dependence becomes more pronounced in network systems with interacting agents and coupling dynamics \cite{nedic2018network}. Each agent selects a local decision, yet its influence can propagate through coupled state dynamics, alter the distributions observed by other agents, and ultimately affect collective performance. 
As such, it is crucial to address these network-mediated effects of local decisions and achieve optimization in an online and distributed manner. 
% It remains unclear how to design a scalable and adaptive algorithm that addresses the network-mediated effects of decisions and achieve closed-loop performance optimization.

% Many real-world decision systems are naturally organized as networks of
% interacting components. In these systems, the network is part of the
% dynamics that generate future states rather than merely a communication
% structure. Local decisions can therefore affect other nodes through the
% coupled dynamics. The resulting steady-state distribution depends on
% both the decisions and the network topology. Optimizing it requires
% accounting for the direct objective derivative and the indirect effect
% through the network response.

% We introduce a centralized network-aware decision-response (NAD)
% gradient and propose distributed NAD (D-NAD), which approximates the
% network-response sensitivity using local and neighboring information.
% We characterize its gradient approximation error and evaluate the method
% on a networked electric-vehicle price--adoption problem.

\subsection{Related Work}
Classical stochastic optimization views the uncertainty distribution as exogenous~\cite{shapiro2009lectures}. In contrast, the recent trend in performative prediction tackles endogenous distribution shifts induced by deployed decisions~\cite{perdomo2020performative,hardt2025performative}.
Typical methods that leverage real-time samples and iterative updates include various stochastic gradient descent algorithms~\cite{mendler2020stochastic,drusvyatskiy2023stochastic}, with extensions to nonconvex~\cite{li2024nonconvex}, bandit-feedback~\cite{chen2024bandit}, and constrained settings~\cite{wang2025constrainedddd,lauand2026online}.

Beyond decision dependency modeled by static parameterized maps, dynamic formulations describe how the distribution evolves over time according to specific dynamics when driven by a (centralized) decision-maker~\cite{cai2024performativecontrol,he2025decisiondependent}.
% Dynamic and performative-control formulations consider distributions
% or states generated by evolving systems and deployed policies.
In multi-agent scenarios, related studies investigate consensus~\cite{li2022multiagent}, network effects~\cite{narang2023multiplayer,wang2023networkeffects}, strategic interactions~\cite{wood2024feedbackgames}, and mechanism design~\cite{gois2025games}. In contrast, we consider cooperative agents with distinct local decisions that jointly determine the steady-state distribution across the network. Hence, our formulation requires capturing the interplay between local decisions and network dynamics, which is not addressed in the aforementioned works.
% which induce network-mediated dynamic distribution shifts and affect a global objective. 

Our work also connects to distributed optimization and networked feedback optimization. In particular, consensus-based distributed optimization coordinates local computations and communication to optimize a global objective over a common decision vector~\cite{yang2019survey}. Moreover, network feedback optimization~\cite{wang2024onlinefeedback} exploits local output measurements and neighborhood communication to optimize the steady-state system operation. Different from these works, we consider a setting where coupled network dynamics generate decision-dependent distributions, and agents select distinct local decisions using samples from the evolving distribution. Further, we leverage local Jacobian approximations as a surrogate for the network-response sensitivity, thereby facilitating an online and distributed solution method.

% and quantify the errors introduced by local network approximations.

% Network Newton further uses neighborhood information and truncated
% matrix expansions to approximate global second-order directions
% ~\cite{mokhtari2017networknewton}. D-NAD adopts a related Neumann
% approximation, but applies it to the network-response sensitivity
% rather than to an inverse Hessian in a consensus formulation.
% \zhmargin{explain differences in problem formulation}

% Networked feedback optimization uses distributed measurements or local
% sensitivities to optimize steady-state operation
% ~\cite{wang2024onlinefeedback}. Unlike
% these methods, we consider decision-dependent distributions and
% quantify errors introduced by local network approximations.

% Recent studies extend these ideas to stochastic,
% price-dependent demand response
% \cite{cavraro2026demandresponse}, while recent surveys
% systematically review the development and stochastic extensions
% of Bass-type diffusion models
% \cite{langrene2026basssurvey}.

\subsection{Contributions}

We study decision-dependent stochastic optimization in networks, where local decisions induce distribution shifts through coupled network dynamics.
% We propose a distributed algorithm through approximating local Jacobians via neighborhood information and a truncated Neumann expansion. 
% We analyze the approximation error of D-NAD and evaluate the
% methods on a networked EV price--demand problem. 
The main contributions are summarized as follows.
\begin{itemize}
    \item We formulate a cooperative stochastic optimization problem, where agents choose distinct local decisions to optimize a global objective under joint decision-dependent steady-state distributions. This formulation captures how local decisions coupled through network interactions and dynamics affect the collective state distribution and the global objective.
    
    % networked decision-dependent stochastic optimization problem, where agents choose local decisions to optimize a global objective under a steady-state distribution induced by coupled network dynamics.
    % problem with coupled network dynamics and introduce the NAD gradient, which captures both direct and network-mediated decision effects.

    \item We propose a distributed network-aware decision-response (D-NAD) gradient algorithm. The key idea is to account for the effect of local decisions on the steady-state network response in iterative adjustments. This is achieved by approximating the network-response Jacobians using local and neighboring information, along with a truncated Neumann expansion. The resulting algorithm is distributed, scalable, and adaptive.
    
    % network-mediated effect of local decisions on the steady-state distribution during iterations. The proposed algorithm features local approximations of the network-response Jacobians using neighborhood information and a truncated Neumann expansion.
    
    \item We establish a non-asymptotic stationarity guarantee and characterize its dependence on the number of iterations and the above approximation errors. We quantify these errors through the Neumann truncation order and the norm bounds of the Jacobians of the network dynamics.
    We demonstrate the performance on a semi-synthetic networked electric-vehicle (EV) price--demand problem calibrated with real-world data.
    
    % Our evaluation is based on these two approximation errors together with the convergence error of centralized NAD (C-NAD). 

    % \item We evaluate D-NAD on a networked EV price--demand problem and study the transient objective error and the effect of the Neumann order. Compared with the performative method, the zero-order gradient method, and local-only NAD (L-NAD) where the nodes only utilize their local information, D-NAD achieves a smaller objective error close to C-NAD.
\end{itemize}

The rest of this paper is organized as follows. Section~\ref{sec:problem_formulation} introduces
the preliminaries and formulates the problem. Section~\ref{sec:distributed_nad} presents
the proposed D-NAD algorithm and its performance guarantees.
Section~\ref{sec:EV_model} reports the numerical results, and Section~\ref{sec:conclusion} concludes this paper. All proofs are provided in the appendix.

\section{Preliminaries and Problem Formulation}
\label{sec:problem_formulation}

We consider a decision-dependent stochastic optimization problem over a network, featuring local distribution shifts driven by collective decisions and coupled through network interactions.
% problem in
% which the decision affects the objective both directly and indirectly
% through the induced distribution of the random state.
Specifically, we consider a network of \(N \in \mathbb{N}_{+}\) interconnected agents.
Each agent \(i\) has a local decision \(u_i \in \mathbb{R}^{n_u}\), a local random state \(x_i \in \mathbb{R}^{n_x}\),
% \zhmargin{call $x_i$ a random variable}
and a local objective function \(\Phi_i(u_i,x_i)\). 
% \zhmargin{dimensionality information}
Define the joint decision and the joint random state as
\(
    u\coloneqq\operatorname{col}(u_1,\ldots,u_N)\in\mathbb{R}^{n_uN}
    \text{ and }
    x\coloneqq\operatorname{col}(x_1,\ldots,x_N)\in\mathbb{R}^{n_xN},
\)
where $n_u$ and $n_x$ are the dimensions of inputs and states respectively.
% here $n$ and $r$ are the dimension of states and inputs at each agent.
% and the feasible set
% \(
%     \mathcal U
%     \coloneqq
%     \mathcal U_1\times\cdots\times\mathcal U_N.
% \)

\subsection{Network Distribution Dynamics}
% For a fixed network matrix \(A\) 
The network-wide joint random state \(x\) evolves according to the following additive stochastic dynamics
\begin{equation}
    x_{k+1}
    =
    f_A(x_k,u_k,d)
    \coloneqq
    \alpha A x_k
    +(1-\alpha)f(x_k,u_k,d),
    \label{eq:network_dynamics}
\end{equation}
where \(k \!\in\! \mathbb{N}\) is the time index, \(A \!\in\! \mathbb{R}^{n_xN\times n_xN}\) is the network interaction matrix, $x_0 \in \mathbb{R}^{n_xN}$ is a random initial state distributed according to $\mu_0$, \(d \!\in\! \mathbb{R}^{n_dN}\) is a stationary random exogenous input or disturbance following the distribution $\mu_d$, 
% \zhmargin{dimension of $d$}
and \(\alpha\!\in\!(0,1)\) balances the linear network coupling and the nonlinear local dynamics.
This additive structure is motivated by classical network models that
combine neighborhood state mixing with local responses
\cite{nedic2018network,proskurnikov2017tutorial}. Here, \(Ax_k\)
captures the state interaction across the network, 
\(f \!\coloneqq\!\operatorname{col}(f_1,\ldots,f_N)\) represents the collection of nonlinear local dynamics of all agents \(\!f_i\!:\! \mathbb{R}^{n_x}\! \times\! \mathbb{R}^{n_u} \!\times \!\mathbb{R}^{n_d} \!\to\! \mathbb{R}^{n_x}\!\), and
% \zhmargin{mind the notation $m$}
% \zhmargin{Does $f$ also have a column structure?}
\(\alpha\) controls the relative influence of these two components.
Such coupled dynamics also arise in practical scenarios including innovation-diffusion and product-adoption models, where local adoption responses interact with social or market spillovers~\cite{kiesling2012agent}.
% ; an EV-adoption instance is
% presented in Sec.~\ref{sec:EV_model}.

We illustrate the interaction structure of the network dynamics~\eqref{eq:network_dynamics} in Fig.~\ref{fig:decentralized_structure}. Here, \(
\mathcal N_i
\!:=\!
\{j\!\neq\! i\!:\!A_{ij}\!\neq\!0\}
\)
denotes the neighborhood of agent \(i\), where \(A_{ij}\!\in\!\mathbb R^{n_x\times n_x}\) is the \((i,j)\)-th block of \(A\). Further, the neighborhood  information set is
\(
\mathcal I_{i,k}^{\mathrm{nei}}
:=
\left\{
\left(
\partial_{x_j}\!f_j,
\partial_{x_j}\!\Phi_j
\right)
\right\}_{j\in\mathcal N_i}
\), where $\Phi_j$ is the objective of agent $j$;  we will shortly specify $\Phi_j$ in Section~\ref{subsec:problem_formulation}. 
The mapping
\(\mathcal K_i\) indicates the local decision rule, which will be given in Section~\ref{sec:distributed_nad}. 
The local component of the exogenous input \(d\) is denoted by \(d_i\).
Moreover, \(f_{A,i}\) indicates the \(i\)-th component of the aggregate network dynamics \(f_A\), and it relies on neighboring states, local decisions and $d_i$. 
The network dynamics cause the shift of the distribution of the local random state \(x_i\) at agent \(i\). As such, agent \(i\) must consider the effects of its own decision \(u_i\) and the decisions of other agents on the resulting distribution shift while collaboratively optimizing the global objective.

% In Fig.~\ref{fig:decentralized_structure},
% The figure is illustrative and introduces no additional dynamics
% beyond \eqref{eq:network_dynamics}.

% To illustrate the local coupling and information structure of
% \eqref{eq:network_dynamics}, Fig.~\ref{fig:decentralized_structure}
% shows a representative agent \(i\), its neighboring agents, and its
% local decision-making structure.

\definecolor{controllerfill}{RGB}{255,241,235}
\definecolor{controllerline}{RGB}{190,65,45}
\definecolor{plantfill}{RGB}{235,244,252}
\definecolor{plantline}{RGB}{45,100,150}
\definecolor{communication}{RGB}{15,125,95}

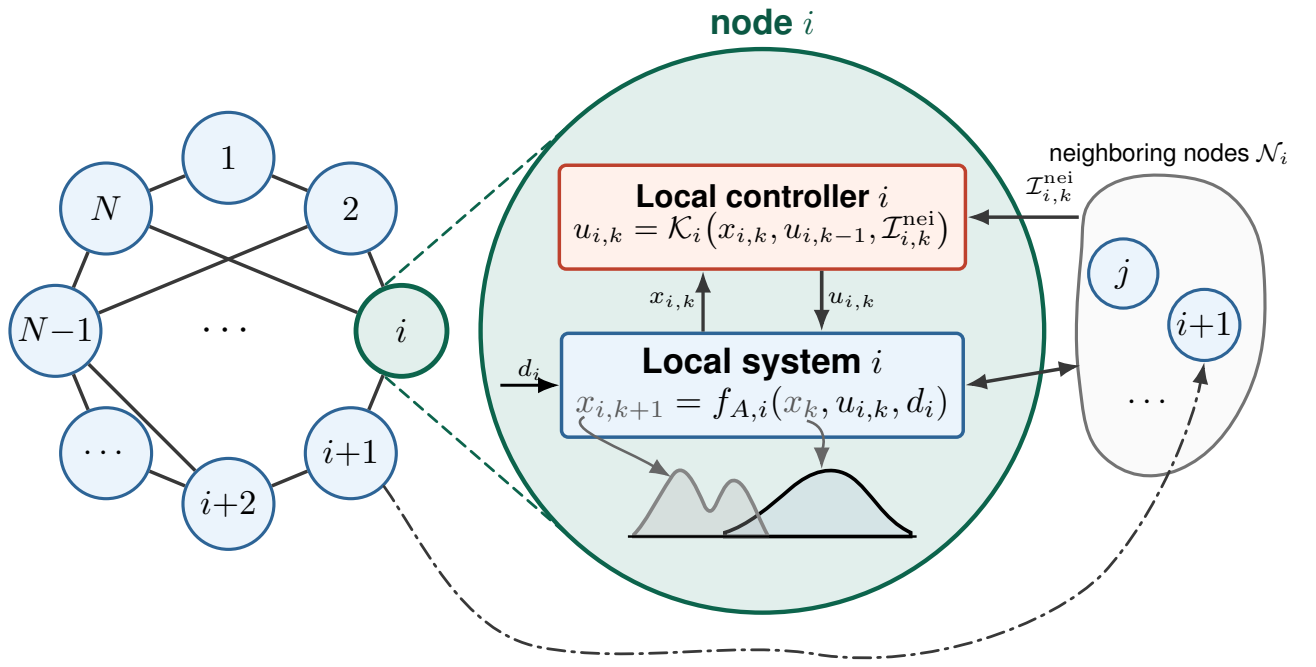
\begin{figure}[!tb]
    \centering
    \setlength{\abovecaptionskip}{-3pt}
    \resizebox{\columnwidth}{!}{%
    \begin{tikzpicture}[
        font=\sffamily\normalsize,
        >={Latex[length=2.2mm]},
        netnode/.style={
            circle, draw=plantline, fill=plantfill,
            line width=0.9pt, minimum size=0.92cm,
            align=center, inner sep=1pt
        },
        selected/.style={
            netnode, draw=communication!80!black,
            fill=communication!12, line width=1.4pt
        },
        controller/.style={
            draw=controllerline, fill=controllerfill,
            rounded corners=2pt, line width=1pt,
            minimum width=3.30cm, minimum height=1.05cm,
            align=center, inner sep=4pt,
            font=\sffamily\small
        },
        dynamics/.style={
            draw=plantline, fill=plantfill,
            rounded corners=2pt, line width=1pt,
            minimum width=3.30cm, minimum height=1.05cm,
            align=center, inner sep=4pt
        },
        networkedge/.style={
            draw=black!78, line width=1pt
        },
        physical/.style={
            <->, draw=black!78, line width=1pt
        },
        information/.style={
            ->, draw=black!78, line width=1pt
        },
        decision/.style={
            ->, draw=black!78, line width=1.05pt
        },
        measurement/.style={
            ->, draw=black!78, line width=0.95pt
        },
        lab/.style={
            font=\sffamily\scriptsize,
            fill=none,
            inner sep=1.2pt
        }
    ]

    % Network overview
    \begin{scope}[shift={(-3.80,0)}]
        \node[netnode] (n1) at (90:1.75cm) {$1$};
        \node[netnode] (n2) at (45:1.75cm) {$2$};
        \node[selected] (ni) at (0:1.75cm) {$i$};
        \node[netnode] (n4) at (-45:1.75cm) {$i{+}1$};
        \node[netnode] (n5) at (-90:1.75cm) {$i{+}2$};
        \node[netnode] (n6) at (-135:1.75cm) {$\cdots$};
        \node[netnode] (n7) at (180:1.75cm) {$N{-}1$};
        \node[netnode] (n8) at (135:1.75cm) {$N$};

        \draw[networkedge] (n1) -- (n2);
        \draw[networkedge] (n2) -- (ni);
        \draw[networkedge] (ni) -- (n4);
        \draw[networkedge] (n4) -- (n5);
        \draw[networkedge] (n5) -- (n6);
        \draw[networkedge] (n6) -- (n7);
        \draw[networkedge] (n7) -- (n8);
        \draw[networkedge] (n8) -- (n1);
        \draw[networkedge] (n2) -- (n7);
        \draw[networkedge] (n8) -- (ni);
        \draw[networkedge] (n5) -- (n7);
        % \draw[networkedge] (n6) -- (n4);

        \node[font=\sffamily\large] at (0,0) {$\cdots$};
    \end{scope}

    % Enlargement guide lines
    \draw[
        draw=communication!80!black,
        line width=0.80pt,
        dash pattern=on 3.2pt off 2.2pt,
        line cap=round
    ]
        (-2.351,0.348)
        --
        (-0.266,2.154);

    \draw[
        draw=communication!80!black,
        line width=0.80pt,
        dash pattern=on 3.2pt off 2.2pt,
        line cap=round
    ]
        (-2.351,-0.348)
        --
        (-0.266,-2.154);

    % Local structure at node i
    \begin{scope}[shift={(0.40,0)}]

        \node[controller] (ctrl) at (1.20,1.15)
        {\textbf{Local controller $i$}\\[-1pt]
        $u_{i,k}
        =
        \mathcal K_i\!\left(
        x_{i,k},
        u_{i,k-1},
        \mathcal I^{\rm nei}_{i,k}
        \right)$};

        \node[dynamics] (dyn) at (1.20,-0.55)
        {\textbf{Local system $i$}\\[-1pt]
        $\textcolor{black!65}{x_{i,k+1}}
        =
        f_{A,i}\!\left(
        \textcolor{black!65}{x_{k}},
        u_{i,k},d_i
        \right)$};

        \node[
            netnode,
            minimum size=0.70cm
        ] (ja) at (4.85,0.58) {$j$};

        \node[
            netnode,
            minimum size=0.70cm
        ] (jb) at (5.67,0.06) {$i{+}1$};

        \node[
            font=\sffamily\normalsize
        ] at (5.12,-0.72) {$\cdots$};

        \node[
            font=\sffamily\scriptsize,
            align=center,
            text=black
        ] at (5.32,1.78)
        {neighboring nodes $\mathcal N_i$};

        % Neighbor information
        \draw[information]
            (4.40,1.15)
            --
            node[
                pos=0.25,
                above,
                yshift=2pt,
                font=\sffamily\scriptsize,
                text=black,
                fill=none,
                inner sep=1.2pt
            ]
            {$\mathcal I^{\rm nei}_{i,k}$}
            (ctrl.east);

        % Physical coupling
        \draw[physical]
            (dyn.east)
            --
            (4.42,-0.35);

        % Decision
        \draw[decision]
            ([xshift=6mm]ctrl.south)
            --
            node[right,lab] {$u_{i,k}$}
            ([xshift=6mm]dyn.north);

        % State measurement
        \draw[measurement]
            ([xshift=-6mm]dyn.north)
            --
            node[left,lab] {$x_{i,k}$}
            ([xshift=-6mm]ctrl.south);

        % Exogenous input
        \draw[
            ->,
            draw=black,
            line width=0.85pt
        ]
            (-1.45,-0.55)
            --
            node[above,lab] {$d_i$}
            (dyn.west);

        % Distribution at iteration k:
        % black solid curve with light-blue shading
        \filldraw[
            draw=black,
            fill=plantline!35,
            fill opacity=0.28,
            draw opacity=1,
            line width=1.05pt,
            line cap=round
        ]
            (0.82,-2.04)
            .. controls (1.05,-1.97) and (1.27,-1.78) ..
            (1.49,-1.59)
            .. controls (1.69,-1.42) and (1.91,-1.35) ..
            (2.10,-1.48)
            .. controls (2.30,-1.62) and (2.42,-1.91) ..
            (2.70,-2.04)
            -- (2.70,-2.08)
            -- (0.82,-2.08)
            -- cycle;

        % Distribution at iteration k+1:
        % gray solid curve with light-gray shading
        \filldraw[
            draw=black!48,
            fill=black!28,
            fill opacity=0.28,
            draw opacity=1,
            line width=1.05pt,
            line cap=round
        ]
            (-0.10,-2.06)
            .. controls (0.02,-1.94) and (0.12,-1.72) ..
            (0.28,-1.48)
            .. controls (0.40,-1.28) and (0.52,-1.55) ..
            (0.62,-1.78)
            .. controls (0.72,-1.94) and (0.75,-1.60) ..
            (0.88,-1.52)
            .. controls (1.02,-1.46) and (1.12,-1.82) ..
            (1.25,-2.06)
            -- (1.25,-2.08)
            -- (-0.10,-2.08)
            -- cycle;

        % Common horizontal axis
        \draw[
            draw=black,
            line width=0.70pt
        ]
            (-0.16,-2.08)
            --
            (2.76,-2.08);

        % From x_{i,k+1} to its distribution
        \draw[
            ->,
            draw=black!60,
            line width=0.75pt
        ]
            ([xshift=-18mm,yshift=-3.2mm]dyn.center)
            .. controls (-0.82,-0.98) and (-0.18,-1.20) ..
            (0.28,-1.48);

        % From x_k to its distribution
        \draw[
            ->,
            draw=black!60,
            line width=0.75pt
        ]
            ([xshift=4.5mm,yshift=-3.6mm]dyn.center)
            .. controls (1.90,-0.97) and (1.80,-1.25) ..
            (1.80,-1.40);

        % Background regions
        \begin{scope}[on background layer]

            \node[
                circle,
                draw=communication!80!black,
                fill=communication!12,
                line width=1.25pt,
                minimum size=5.70cm,
                label={
                    [font=\sffamily\bfseries,
                    text=communication!70!black]
                    90:node $i$
                }
            ] at (1.20,0) {};

            \filldraw[
                fill=black!2,
                draw=black!55,
                line width=0.9pt
            ]
                plot[
                    smooth cycle,
                    tension=0.72
                ]
                coordinates {
                    (4.42,0.72)
                    (4.58,1.38)
                    (5.32,1.57)
                    (6.03,1.27)
                    (6.28,0.44)
                    (6.10,-0.66)
                    (5.42,-1.40)
                    (4.66,-1.28)
                    (4.40,-0.46)
                };

        \end{scope}
    \end{scope}

    % Correspondence between the two i+1 nodes
    \draw[
        ->,
        draw=black!78,
        line width=0.80pt,
        dash pattern=on 5pt off 2pt on 0.8pt off 2pt,
        line cap=round
    ]
        (n4.south east)
        .. controls (-1.85,-2.10) and (-1.45,-3.05) ..
        (-0.65,-3.20)
        .. controls (0.30,-3.38) and (1.35,-3.12) ..
        (2.20,-3.28)
        .. controls (3.25,-3.42) and (4.50,-3.00) ..
        (5.05,-2.45)
        .. controls (5.60,-1.85) and (5.95,-0.95) ..
        (jb.south);

    \end{tikzpicture}
    }
    \caption{This figure illustrates the interaction structure of the network system. To achieve collective performance optimization, each agent \(i\) adjusts its decision by using local samples, computation, and communication with its neighbors and addressing the effect of decision dependence.
    % interacts with neighboring agents
    % through the network dynamics and uses locally available neighbor
    % information to determine its decision.
    }
    \label{fig:decentralized_structure}
\end{figure}
% \zhmargin{visualize distribution shifts}

\subsection{Problem Formulation}\label{subsec:problem_formulation}
We impose the following technical assumptions and formally introduce the problem of interest.

% We impose the following assumptions on the network dynamics, the
% steady-state objective, and the stochastic gradient.

\begin{assumption}
\label{asm:network-dynamics}
The network is undirected, and the matrix \(A\) is nonnegative and
doubly stochastic. The dynamics \(\!f(x,u,d)\!\) are continuously
differentiable, \(L_f^x\)-Lipschitz continuous in \(x\), \(L_f^u\)-Lipschitz continuous in \(u\), and
\(L_f^d\)-Lipschitz continuous in \(d\), where
\(L_f^x\!\in\!(0,1)\), \(L_f^u\!>\!0\), and \(L_f^d\!>\!0\). There exists a unique and continuously differentiable steady-state
response map \(h_A\) such that $h_A\!(u,\!d)\!=\!f_A\!\bigl(h_A\!(u,d),\!u,\!d\bigr)$.
% The steady-state
% response map \(h_A\) defined in \eqref{eq:steady_state_map} is .
Furthermore, for all \(u\), \(\nabla_x f(x,u,d)\) and
\(\nabla_u f(x,u,d)\) are \(M_f^x\)-Lipschitz and
\(M_f^u\)-Lipschitz continuous in \((x,d)\), respectively.
\end{assumption}

% For a fixed decision \(u\) and a realization of \(d\), Assumption~1
% guarantees a unique steady-state response \(h_A(u,d)\) satisfying
% \begin{equation}
%     h_A(u,d)
%     =
%     f_A\bigl(h_A(u,d),u,d\bigr).
%     \label{eq:steady_state_map}
% \end{equation}

Assumption~\ref{asm:network-dynamics} adapts the contraction and regularity conditions in
\cite{he2025decisiondependent} to the setting of coupled network dynamics~\eqref{eq:network_dynamics}. Contraction ensures
that the state admits a unique steady-state response for
fixed \(u\) and \(d\), while the differentiability and Jacobian
regularity support the sensitivity and Neumann-series approximations to be presented.
The nonnegative doubly stochastic matrix and undirected-network conditions
are standard in network averaging models
\cite{proskurnikov2017tutorial,nedic2018network}.

Given a fixed joint decision \(u\), the local steady state \(x_{i,\mathrm{ss}}=h_{A,i}(u,d)\) follows the distribution
\(
\mu_i^{\mathrm{ss}}(u)
:=h_{A,i}(u,\cdot)_{\#}\mu_d,
\)
i.e., the pushforward of \(\mu_d\) under \(h_{A,i}(u,\cdot)\). 
Without loss of generality, let
\(
\mathbb E_{(x_0,d)}
[\|x_0-h_A(0,d)\|]<\infty
\) hold,
where the expectation is taken over the joint distribution of
 \(x_0\) and \(d\).
% The distribution \(\mu_i(u)\) in
% \eqref{eq:network_ddso} is the distribution of
% \(
% x_{i,\mathrm{ss}}
% =
% h_{A,i}(u,d)
% \)
% induced by \(d\sim\mu_d\).
The network-wide performance optimization problem is
\begin{equation}
    \mathop{\mathrm{min}}_{u \in \mathbb{R}^{n_uN}}
    \;
    \sum_{i=1}^{N}
    \mathbb E_{x_i \sim \mu_i^{\text{ss}}(u)}
    \left[\Phi_i(u_i,x_i)\right].
    \label{eq:network_ddso}
\end{equation}
% where \(\mu_i^\text{ss}(u)\) is the steady-state distribution of the local random state
% \(x_i\) . 
Although the local decision \(u_i\) directly enters only the local dynamics of agent \(i\), network coupling causes it to affect the states and objectives of other agents. Moreover, thanks to the contraction property in Assumption~\ref{asm:network-dynamics}, the steady-state distribution \(\mu_i^{\mathrm{ss}}(u)\) is not affected by the initial distribution. Consequently, \(\mu_i^{\mathrm{ss}}(u)\) is determined by the joint decision \(u\) and the input law \(\mu_d\). Using the steady-state response map \(h_A\), with \(\Phi(u,x):=\sum_i\Phi_i(u_i,x_i)\), the objective in
\eqref{eq:network_ddso} can be equivalently written as
\begin{equation}
\widetilde{\Phi}(u)
:=
\mathbb E_{d\sim\mu_d}
\left[
\Phi\bigl(u,h_A(u,d)\bigr)
\right].
\label{eq:steady_state_objective}
\end{equation}
% which is the steady-state objective.

\begin{assumption}
\label{asm:objective-regularity}
The steady-state objective \(\widetilde{\Phi}\) in \eqref{eq:steady_state_objective}
is well-defined, \(L_{\widetilde{\Phi}}\)-smooth
(with \(L_{\widetilde{\Phi}}\)-Lipschitz gradients), and bounded below.
The function \(\Phi(u,x)\) is \(L_\Phi^x\)-Lipschitz continuous in
\(x\). The partial gradients \(\nabla_u\Phi(u,x)\) and
\(\nabla_x\Phi(u,x)\) are \(M_\Phi^u\)-Lipschitz and
\(M_\Phi^x\)-Lipschitz in \(x\), respectively.
\end{assumption}

Assumption~\ref{asm:objective-regularity} follows the regularity conditions used in
\cite{he2025decisiondependent}. These conditions are standard and satisfied in practical applications~\cite{shapiro2009lectures,bottou2018optimization}. A well-defined, smooth, and
lower-bounded objective facilitates convergence analysis, while the
Lipschitz conditions constrain the gradient error caused by using
transient rather than steady-state responses. 
% \zhmargin{Support assumptions not by explaining how they are used, but relating them to practical problems and references}

\begin{assumption}
\label{asm:lipschitz-domination}
There exists a positive random variable \(L_{\Phi}^u\) such that
\(
\mathbb E_{d\sim\mu_d}[L_{\Phi}^u(d)]\!<\!\infty,
\)
and, for all \(u_1,u_2\!\in\!\mathbb R^{n_uN}\),
\[
\left|
\Phi(u_1,h_A(u_1,d))
\!-\!
\Phi(u_2,h_A(u_2,d))
\right|
\leq
L_{\Phi}^u(d)\|u_1\!-\!u_2\|.
\]
\end{assumption}

Assumption~\ref{asm:lipschitz-domination} is
standard in stochastic programming
\cite{shapiro2009lectures}. The integrable random Lipschitz bound
provides the domination required to interchange the expectation and
gradient operators when differentiating the steady-state objective.

To collectively solve~\eqref{eq:network_ddso}, each agent must capture both the direct
effect of the local decision and the indirect effect propagated
through the network. Since no individual agent could possess complete network information, each agent must leverage local computation and communication. This motivates our distributed network-aware algorithm.

\section{Distributed Network-Aware Decision-Response (D-NAD) Gradient Algorithm}
\label{sec:distributed_nad}

We propose a distributed network-aware decision-response gradient algorithm. To provide intuition, we first present a centralized version that exploits complete network information. 
% For each sample \(m\) at iteration \(k\), define the gradients
% \(H_k^m:=\nabla_x f_A \in \mathbb{R}^{n_xN\times n_xN}\),
% \(B_k^m:=\nabla_u f_A \in \mathbb{R}^{n_xN\times n_uN}\),
% \(p_k^m:=\nabla_u\Phi \in \mathbb{R}^{n_uN}\), and
% \(q_k^m:=\nabla_x\Phi \in \mathbb{R}^{n_xN}\) evaluated at the current transient state $x_{k}^m$, decision $u_{k}$, and exogenous input $d^m$. 
Following the gradient construction in~\cite[Section~3]{he2025decisiondependent},
we define the centralized NAD (C-NAD) update
\begin{subequations}
\label{eq:centralized_nad_gradient}
\begin{align}
\!&u_{k+1}
\!=\!
u_k\!-\!\eta\,\widehat g_k^{\mathrm{C\text{-}NAD}}\!,\quad \!
\widehat g_k^{\mathrm{C\text{-}NAD}}
\!\!=\!
\tfrac{1}{N_{\mathrm{mb}}}
{\textstyle\sum\nolimits_{m=1}^{N_{\mathrm{mb}}}}
\widehat g_{k,m}^{\mathrm{C\text{-}NAD}}\!,
\label{eq:centralized_gd_a}
\\
&\begin{aligned}[c]
\widehat g_{k,m}^{\mathrm{C\text{-}NAD}}
&\!=
\!\nabla\!_u\Phi(u_k,x_k^m)
\!+\!
\nabla\!_u f\!_A(x_k^m,u_k,d^m)^\top\\
&\;\;\quad\cdot
(I\!-\!\nabla\!_x f\!_A(x_k^m,u_k,d^m))\!^{-\!\top}
\!\nabla\!_x\Phi(u_k,x_k^m).
\end{aligned}
\label{eq:centralized_gd}
\end{align}
\end{subequations}
\noindent
In \eqref{eq:centralized_gd}, $N_{\mathrm{mb}}$ is the mini-batch size, and $(x_k^m, d^m)$ represent the $m$-th sample pair of the random state at time $k$ and the exogenous input. Each sample pair may correspond to an individual record in a dataset or a realized outcome from a simulator run and is propagated through the network dynamics. 
Further, \(\nabla_u\Phi\) represents the direct effect of the decision on the objective, while the second term captures its indirect effect through the steady-state network response. The steady-state sensitivity \(\nabla_u h_A := (I-\nabla_x f_A)^{-1}\nabla_u f_A\) follows from implicit differentiation of the steady-state equation~\cite[Theorem~9.28]{rudin1976principles}, with the derivatives evaluated at \(x=h_A(u,d)\). Therefore, \eqref{eq:centralized_nad_gradient} accounts for both the immediate objective variation and the decision response propagated through the coupled network dynamics.

However, the C-NAD method cannot be directly executed by individual
agents in the distributed setting. In particular, computing
\((I\!-\!\nabla_x f\!_A)^{-1}\) requires the complete network Jacobian \(\nabla_x f\!_A\),
which depends on the global network topology and the dynamic derivatives
of all agents, as well as the partial gradient vector \(\nabla_x \Phi\).
Nonetheless, an individual agent \(i\) only has access to its own local derivatives
and the information received from its neighbors. Hence, \eqref{eq:centralized_nad_gradient} cannot be evaluated exactly from the prescribed local information.

\subsection{D-NAD Algorithm}

\setcounter{algocf}{0}
\begin{algorithm}[t]
\caption{Distributed NAD Gradient Algorithm}
\label{alg:distributed-neumann-nad}

\KwIn{
Neumann order $R$, mini-batch size $N_{\mathrm{mb}}$,
initial decisions $\{u_{i,0}\!=\!0^{n_u}\}$, step sizes $\{\eta_i\}$.
}

\For{iteration $k=0,1,2,\ldots$}{

    Evolve the network distribution dynamics:\label{line:propagate}
    % $m\in\{1,\ldots,N_{mb}\}$:
    \(
    \begin{aligned}[t]
    \mu_{k+1}
    \!=\!
    f_{A}(\cdot,u_k,d)_{\#}\mu_{k},
    \end{aligned}
    \)
    % m=1,\ldots,N_{\mathrm{mb}}

    Draw $N_{\mathrm{mb}}$ samples from the distributions: \(
    \begin{aligned}[t]
        x_k^m\sim\mu_k,\quad m=1,\ldots,N_{\mathrm{mb}}
    \end{aligned}
    \)\label{line:draw_mb_samples}

    \For{agent $i=1,\ldots,N$}{
        Compute dynamics and objective derivatives:\newline
        \(
        c_{i,k}^m
        \!=\!
        \frac{\partial f_i}{\partial x_i},
        b_{i,k}^m
        \!=\!
        \frac{\partial f_i}{\partial u_i},
        p_{i,k}^m
        \!=\!
        \frac{\partial\Phi_i}{\partial u_i},
        q_{i,k}^m
        \!=\!
        \frac{\partial\Phi_i}{\partial x_i}.
        \)\label{line:derivatives}
    
        Receive
        \(
        \{c_{j,k}^m,q_{j,k}^m\}_{j\in\mathcal N_i}
        \)
        from its neighbors. \label{line:communicate}

        Construct the neighborhood Jacobian:
        \(
        \begin{aligned}[t]
        H_{\mathcal S_i,k}^m
        ={}
        \alpha A_{\mathcal S_i}
        +
        (1-\alpha)\operatorname{blkdiag}
        \bigl(\{c_{j,k}^m\}_{j\in\mathcal S_i}\bigr).
        \end{aligned}
        \)\label{line:jacobian}

        Initialize:
        \(
        v^m
        =
        y^m
        =
        \operatorname{col}(\{q_{j,k}^m\}_{j\in\mathcal S_i}).
        \)
        
        \For{$r=1,\ldots,R$\label{line:neumann-loop}}{
            \(
            v^m
            =
            \big({H_{{\mathcal S_i,k}}^m}\big)^\top v^m,
            \)\label{line:neumann-v}
            
            \(
            y^m
            =
            y^m
            +
            v^m.
            \)\label{line:neumann-w}
        }

        Compute the distributed NAD gradient: \label{line:nad-gradient}
        \(
        \begin{aligned}[t]
        \widehat g_{i,k}^{\textrm{D-NAD}}
        \!=\!
        \tfrac{1}{N_{\mathrm{mb}}}
        {\textstyle\sum\nolimits_{m=1}^{N_{\mathrm{mb}}}}
        \left[
        p_{i,k}^m
        +(1-\alpha){b_{i,k}^m}^\top
        y_i^m
        \right].
        \end{aligned}
        \)

        Update:
        \(
        u_{i,k+1}
        =
        u_{i,k}
        -
        \eta_i\,\widehat g_{i,k}^{\textrm{D-NAD}}.
        \) \label{line:nad-update}
    }
}
\end{algorithm}

To overcome the above challenges, we
propose the distributed network-aware decision-response (D-NAD) algorithm. The main
idea is to approximate the unavailable global network-response sensitivity using
local and neighboring Jacobians.
Specifically, each
agent constructs a local representation of the network
Jacobian through communication with its neighbors and uses a truncated Neumann expansion to avoid directly
inverting the resulting matrix. 

Algorithm~\ref{alg:distributed-neumann-nad} summarizes the complete D-NAD updates. In line~\ref{line:propagate}, the distribution dynamics are propagated through the network, keeping exogenous input \(d\) fixed over time.
Line~\ref{line:draw_mb_samples} draws $N_\text{mb}$ independent samples from the propagated distributions.
Line~\ref{line:derivatives} involves the computation of the local partial derivatives,
\(
c_{i,k}^m\!:=\!\partial_{x_i} f_i\!\in\!\mathbb{R}^{n_x\times n_x},
\)
\(
b_{i,k}^m\!:=\!\partial_{u_i} f_i\!\in\!\mathbb{R}^{n_x\times n_u},
\)
\(
p_{i,k}^m\!:=\!\partial_{u_i}\Phi_i\!\in\!\mathbb{R}^{n_u},
\)
\(
q_{i,k}^m\!:=\!\partial_{x_i}\Phi_i\!\in\!\mathbb{R}^{n_x}
\) evaluated at $(x_{i,k}^m, u_{i,k}, d_i^m)$ analogously to~\eqref{eq:centralized_nad_gradient}. 
In line~\ref{line:communicate}, local derivatives are communicated across agents.
In line~\ref{line:jacobian}, $\mathcal N_i$ denotes the neighborhood of agent $i$, and \(\mathcal S_i\!:=\!\{i\}\!\cup\!\mathcal N_i\) denotes its closed
neighborhood. The matrix
\(A_{\mathcal S_i}\!\in\!\mathbb R^{n_x|\mathcal S_i|\times n_x|\mathcal S_i|}\)
is obtained from \(A\)
by retaining the rows and columns corresponding to agent \(i\) and its direct neighbors and adding the omitted row weights to the corresponding diagonal entries. Unlike distributed optimization methods~\cite{nedic2018network,yang2019survey} that perform averaging across the network to reach consensus, the proposed algorithm leverages the submatrix $A_{\mathcal S_i}$ to construct a local Jacobian approximation, thereby facilitating local updates. 
Specifically, lines~\ref{line:neumann-loop}--\ref{line:neumann-w} represent the inner loop of approximating the network-response sensitivity $(I-H_{\mathcal S_i,k}^m)^{-\top}
q_{\mathcal S_i,k}^m$ via a truncated Neumann expansion $\sum_{r=0}^{R}
\big({H_{\mathcal S_i,k}^m}^\top\big)^r
q_{\mathcal S_i,k}^m$.
For agent $i$, the complexity of constructing this approximation is $\mathcal O(R|\mathcal S_i|^2)$, whereas the complexity of computing the original sensitivity is $\mathcal O(|\mathcal S_i|^3)$. Hence, by selecting a suitable truncation order \(R\), the computational cost can be significantly reduced; see Section~\ref{sec:EV_model} for details. Finally, in line~\ref{line:nad-gradient}, the local mini-batch gradient is computed and used to update the decision via line~\ref{line:nad-update}.

Algorithm~\ref{alg:distributed-neumann-nad} features local iterative updates using samples from the evolving network system and communication of local derivatives. Furthermore, the D-NAD algorithm can be extended to i) handle input constraints by incorporating a projection step in line~\ref{line:nad-update}, and ii) address non-Euclidean geometries via a mirror descent framework.

% In this algorithm, each agent updates only its own decision variable. Neighbor
% communication is used to exchange local derivatives rather than to
% average the decisions, and therefore no decision consensus is required.
% The derivatives are exchanged once per optimization iteration, while
% \(R\) controls only the amount of local computation.

% \zhmargin{briefly mention that other extensions (e.g., projected descent) are available.}

\subsection{Approximation Error and Convergence Analysis}
Compared to its centralized benchmark~\eqref{eq:centralized_nad_gradient}, the D-NAD algorithm leverages an approximate neighborhood Jacobian and a truncated Neumann expansion; see line~\ref{line:jacobian} and lines \ref{line:neumann-loop}--\ref{line:neumann-w} of Algorithm~\ref{alg:distributed-neumann-nad}, respectively. 
We quantify the effects of these approximations on the distributed NAD gradient through the Neumann truncation error and the neighborhood approximation error. Then, we combine these error bounds to characterize the convergence properties of D-NAD\@. 

% Relative to C-NAD, the D-NAD
% direction introduces a Neumann truncation error and a neighborhood
% approximation error. We first derive separate bounds for these two
% components and then combine them to bound the total D-NAD gradient
% error. Finally, this error bound is incorporated into the convergence
% result of C-NAD to establish a finite-time stationarity bound for
% D-NAD.

For every \(\gamma\in[0,1)\), we define
\begin{subequations}
\label{eq:kappa}
\begin{align}
&\kappa_R(\gamma)
:=
{\textstyle\sum\nolimits_{r=2}^{R}}r\gamma^{r-1}
=
\tfrac{1-(R+1)\gamma^R+R\gamma^{R+1}}{(1-\gamma)^2}-1,
\label{eq:kappa_R}
\\
&\begin{aligned}[c]
\kappa_\infty\!(\gamma)
:=
{\textstyle\sum\nolimits_{r=2}^{\infty}}r\gamma^{r-1}
=
\tfrac{1}{(1-\gamma)^2}-1.
\end{aligned}
\label{eq:kappa_inf}
\end{align}
\end{subequations}
Let
\(
\Delta_A:=\sum_{i=1}^{N}\|\Delta_i\|^2
\)
denote the network submatrix approximation error, where
\(
\Delta_i:=A-\widetilde A_{\mathcal S_i}.
\)
Here, $\widetilde A_{\mathcal S_i}$ is the lifting of \(A_{\mathcal S_i}\) to the full network dimension, obtained by filling all off-diagonal entries of non-neighboring agents with zeros and their diagonal entries with ones. 
Analogously,
\begin{equation}
\widetilde H_{\mathcal S_i,k}^m
\!:=\!
\alpha\widetilde A_{\mathcal S_i}
\!+\!(1\!-\!\alpha)
\operatorname{blkdiag}(c_{1,k}^m,\ldots,c_{N,k}^m)
\label{eq:extended_H}
\end{equation}
is the lifting of \(H_{\mathcal S_i,k}^m\).
Its norm is bounded above by
\(\|\widetilde H_{\mathcal S_i,k}^m\|\leq\alpha+(1-\alpha)L_f^x<1\)
under Assumption~\ref{asm:network-dynamics}.
Define a common Jacobian norm upper bound for each iteration \(k\)
\begin{equation}
\gamma_k
:=\!\!
\max_{1\leq m\leq N_{\mathrm{mb}}}\!\!
\left\{
\|H_k^m\|,
\max_{1\leq i\leq N}
\|\widetilde H_{\mathcal S_i,k}^m\|
\right\},
\label{eq:gamma_k_definition}
\end{equation}
% for each iteration \(k \in \mathbb{N}\).
where \(H_k^m\) denotes \(\nabla\!_x f\!_A(x_k^m,u_k,d^m)\). 
By Assumption~1 and \eqref{eq:gamma_k_definition}, the Jacobians are contractive, and hence
\(
0\leq\gamma_k<1.
\)

With \(B_k^m\) denoting \(\nabla\!_uf\!_A(x_k^m,u_k,d^m)\),
% For each iteration $k \in \mathbb{N}$, 
we define \[
e_k^{(1)}(R)
:=
\tfrac{1}{N_{\mathrm{mb}}}
{\textstyle\sum\nolimits_{m=1}^{N_{\mathrm{mb}}}
(B_k^m)^\top
\sum_{r=R+1}^{\infty}
((H_k^m)^\top)^r q_k^m}
\]
to be the transient Neumann truncation error. 
Similarly, the local transient neighborhood approximation error is
\[
e_{i,k}^{(2)}(\!R)
\!:=\!
\tfrac{1}{N_{\mathrm{mb}}}
{\textstyle\sum\nolimits_{m=1}^{N_{\mathrm{mb}}}
\!(\!B_{ii,k}^m)\!^\top
\!\!\sum_{r=2}^{R}
\!\left[
\!\left(\!\!
(\!H_k^m)^{r}
\!\!\!\!-\!
(\!\widetilde H_{\mathcal S_i,k}^m)^{r}
\!\right)\!\!^\top \!\!q_k^m
\right]_i},
\]
where \(B_{ii,k}^m\) is the \((i,i)\)-th block of \(B_{k}^m\).
The transient neighborhood error is
\(
e_k^{(2)}\!(R)
\!=\!
\operatorname{col}\bigl(e_{1,k}^{(2)}(R),\ldots,e_{N,k}^{(2)}(R)\bigr).
\)
The sum starts at \(r\!=\!2\) because the terms for \(r\!=\!0\)
and \(r\!=\!1\) vanish. The former uses identical identity
matrices, while the latter follows from the identical
\(i\)-th block columns of \(H_k^m\) and
\(\widetilde H_{\mathcal S_i,k}^m\), since
\(\mathcal S_i\) contains agent \(i\) and all its neighbors.
For notational simplicity we define
\(
C_0:=(1-\alpha)^2(L_f^uL_\Phi^x)^2.
\)

% The total D-NAD gradient error is
% \(
% e_{k}^{\mathrm{D\text{-}NAD},R}
% :=
% e_{k}^{k,R}+e_{\mathrm{nei}}^{k,R}.
% \)

% For index
% \(\mathord{\scriptstyle\times}\!\in\!
% \{(1),\!(2),\!\mathrm{C\text{-}NAD},\!\mathrm{D\text{-}NAD}\}\),
% define the errors
% \(\mathcal E^{\times}(T,R):=
% \frac1T\sum_{k=0}^{T-1}\mathbb E[\|e_k^{\times}\|^2]\),
% where \(\mathcal E^{(1)}\) and \(\mathcal E^{(2)}\) stand for the
% Neumann truncation error and the neighborhood approximation error,
% respectively. For notational convenience, we write
% \(\mathcal E^{\times}\) for \(\mathcal E^{\times}(T,R)\) hereafter.

\begin{theorem}[Neumann truncation error]
\label{lem:average-neumann-error}

Let Assumptions~\ref{asm:network-dynamics} and~\ref{asm:objective-regularity} hold. For any $T\geq1,R\geq1$,
% \(\forall T\!\geq\!1\) and \(\forall R\!\geq\!1\),
the Neumann truncation error satisfies
\begingroup
\[
\frac1T\sum_{k=0}^{T-1}\mathbb E[\|e_k^{(1)}\|^2]
\leq
C_0\mathcal B^{(1)},
\]
\endgroup
where 
\(
\!\mathcal B^{(1)}
\!:=\!
\tfrac1T\!
\sum_{k=0}^{T-1}
\!\mathbb E\!
\left[
\tfrac{\gamma_k^{2R+2}}
{(1-\gamma_k)^2}
\right]\!
\)
decreases with $R$.
\end{theorem}

Theorem~\ref{lem:average-neumann-error} implies that \(
\lim\nolimits_{R\to\infty}\!
\tfrac{\gamma_k^{2R+2}}{(1-\gamma_k)^2}
\!=\!0,
\)
$\forall\gamma_k\!\in\! [0,1)$.
Therefore, $\mathcal B^{(1)}$ approaches zero as $R$ increases to infinity. In other words, the Neumann truncation error can be made arbitrarily small by increasing
\(R\). A smaller \(\gamma_k\) allows a lower $R$ to achieve a prescribed truncation accuracy.

\begin{theorem}[Neighborhood approximation error]
\label{lem:average-neighborhood-error}

Let Assumptions~\ref{asm:network-dynamics} and~\ref{asm:objective-regularity} hold. For any $T\geq1, R\geq1$,
% \(\forall T\geq1\) and \(\forall R\geq1\),
the neighborhood approximation error satisfies
\[
\frac1T\sum_{k=0}^{T-1}\mathbb E[\|e_k^{(2)}\|^2]
\leq
C_0\mathcal B^{(2)}
\leq
C_0\mathcal B^{(2)}_\infty,
\]
where
\(
\mathcal B^{(2)}
\!:=\!
\frac{\alpha^2\Delta_A}{T}
\!\!\sum_{k=0}^{T-1}
\mathbb E[\kappa_R(\gamma_k)^2]
\)
increases with $R$ and converges to
\(\mathcal B^{(2)}_\infty\), which
is defined analogously by replacing
\(\kappa_R(\gamma_k)\) with \(\kappa_\infty(\gamma_k)\) according to \eqref{eq:kappa}.
\end{theorem}

Theorem~\ref{lem:average-neighborhood-error} indicates that, even though \(\mathcal B^{(2)}\) grows as the Neumann order $R$ increases, it remains bounded by \(\mathcal B^{(2)}_\infty\). Therefore, the neighborhood approximation error remains bounded. For a fully connected network, \(\Delta_A\!=\!0\) and the neighborhood error vanishes. Overall, the above bounds on \(\mathcal B^{(1)}\) and \(\mathcal B^{(2)}\) reveal a trade-off in selecting the Neumann order \(R\): a small \(R\) saves computational effort but may leave a large truncation error. Increasing \(R\) reduces this error, while the neighborhood error bound grows toward a finite limit. Therefore, a suitable \(R\) should balance approximation accuracy and computational cost. We will empirically evaluate the performance of D-NAD with various $R$ in Section~\ref{sec:EV_model}. 
% as an excessively large \(R\) increases computational cost without necessarily improving the overall accuracy.

% \zhmargin{The definitions of $e_k$ are scattered}
% For both C-NAD and D-NAD, the true gradient at iteration $k$ can be
% decomposed into a gradient estimate and its error, i.e.,
% \(
% \nabla\widetilde{\Phi}(u_k)
% =
% \widehat g_k^{\mathrm{C\text{-}NAD}}
% +
% e_k^{\mathrm{C\text{-}NAD}}
% =
% \widehat g_k^{\mathrm{D\text{-}NAD}}
% +
% e_k^{\mathrm{D\text{-}NAD}}.
% \)
% To facilitate convergence analysis, we make the following assumption on the stochastic gradient variance.
\begin{assumption}
\label{asm:gradient-variance}
The stochastic gradient satisfies
\begingroup
\setlength{\abovedisplayskip}{9pt}
\setlength{\abovedisplayshortskip}{9pt}
\setlength{\belowdisplayskip}{9pt}
\setlength{\belowdisplayshortskip}{9pt}
\[
\begin{aligned}
\operatorname{Var}\!
\left[
\widehat g_{k,m}^{\mathrm{C\text{-}NAD}}
\,\!\middle|\!\,
\mathcal F_{k-1}
\right]
\!\leq\!
M
\!+\!\!
M_V
\!\left\|
\bar g_k^{\mathrm{C\text{-}NAD}}
\right\|^2\!\!,\;
m\!=\!1,\!...,\!N_{\mathrm{mb}},
\end{aligned}
\]
\endgroup
where \(M,M_V\geq0\), \(N_{\mathrm{mb}} \in \mathbb{N}_{+}\) is the mini-batch size, and \(\mathcal F_{k-1}\) denotes all information available before iteration \(k\).
\end{assumption}
Here,
\(
\bar g_k^{\mathrm{C\text{-}NAD}}
:=
\mathbb E[
\widehat g_k^{\mathrm{C\text{-}NAD}}
\mid\mathcal F_{k-1}]
\)
is the conditional expectation of the stochastic C-NAD gradient defined
in~\eqref{eq:centralized_gd_a}. Since its network-response sensitivities are
evaluated at the current transient states, it generally differs from
the true steady-state gradient \(\nabla\widetilde{\Phi}(u_k)\).
% ; see also \eqref{eq:centralized_gd} in the following section.
% \zhmargin{discuss this notation}
Assumption~\ref{asm:gradient-variance} uses the relative variance condition in
\cite{he2025decisiondependent}, which is also standard in stochastic
gradient analysis \cite{bottou2018optimization}. It allows nonzero
gradient variance at a stationary point while restricting its growth, thereby facilitating the following convergence characterization.

\begin{theorem}[Convergence characterization of D-NAD]
\label{thm:distributed-nad-convergence}

Suppose that Assumptions~\ref{asm:network-dynamics}--\ref{asm:gradient-variance}
hold.
Choose a common step size
\(
\eta_i\!=\!\eta\!=\!\tfrac{\eta_0}{\sqrt{T}},
\)
with \(\eta_0\!>\!0\) sufficiently small.
For every \(T\!\geq\!1\) and \(R\!\geq\!1\), the D-NAD iterates \(\{u_k\}_{k=0}^{T-1}\) satisfy
\begin{align}
    \frac1T
    \sum_{k=0}^{T-1}
    \mathbb E
    \left[
    \left\|
    \nabla\widetilde{\Phi}(u_k)
    \right\|^2
    \right]
    =
    \mathcal O\!\left(\frac1{\sqrt T}+
    \mathcal B^{(1)}
    +
    \mathcal B^{(2)}
    \right).
\label{eq:avg_2nd_moment}
\end{align}
\end{theorem}

Theorem~\ref{thm:distributed-nad-convergence} establishes a non-asymptotic bound on the convergence property of D-NAD and characterizes its dependence on the number of iterations \(T\) and the Neumann order \(R\). The scaling of the first term in the bound with respect to \(T\) matches that of stochastic gradient methods for nonconvex optimization~\cite{bottou2018optimization}, here extended to a setting with decision-dependent network dynamics. 
Increasing \(R\) reduces the Neumann truncation error bound \(\mathcal B^{(1)}\), while the neighborhood approximation error bound \(\mathcal B^{(2)}\) accumulates but remains bounded by \(\mathcal B^{(2)}_\infty\), implying that the overall convergence effect of the two approximations remains controlled. In practice, we typically choose \(R\!=\!2\) or \(3\) to balance approximation accuracy and computational cost; see also Section~\ref{sec:EV_model}. For fixed \(R\) and neighborhoods, increasing \(T\) reduces the iteration-dependent term, while the approximation terms determine the remaining stationarity bound.

% \zh{Add discussions of the significance and implications of Theorem~\ref{thm:distributed-nad-convergence}.}

\section{EV Price-Demand Case Study}
\label{sec:EV_model}

We consider a networked EV price--demand model with $N$ agents, where each agent
represents the market condition of a country. The state
\(
x_{i,k}\in[0,1]^3
\)
denotes the EV adoption level, while
\(
u_{i,k}\in \mathbb{R}_{+}^{3}
\)
is the local price decision of agent \(i\) at time \(k\). The local price decision affects the market adoption level mainly through an exponential demand response mapping \(\psi_i(u_i)=a_i\exp(-b_i u_i)\), where \(a_i > 0\) and \(b_i > 0\) represent the demand scale and the price sensitivity, respectively. 
In the experiments, we consider \(N=50\) agents, each with three states corresponding to battery electric vehicle, plug-in hybrid electric vehicle, and fuel cell electric vehicle. 
% \zhmargin{These abbreviations are unclear}
We construct the semi-synthetic dataset by replicating selected country-level EV sales and price data\footnote{See 
\href{https://www.iea.org/data-and-statistics/data-tools/global-ev-data-explorer}
{IEA Global EV Data: https://www.iea.org/data-and-statistics/data-tools/global-ev-data-explorer} and  
\href{https://ourworldindata.org/grapher/price-of-lithium-ion-battery-cells}
{Our World in Data lithium-ion battery-cell prices: https://ourworldindata.org/grapher/price-of-lithium-ion-battery-cells}.}
across the agents and rescaling their fitted parameters.

% \zhmargin{These links will be removed in the submission system. Make them appear as urls.}

% \zhmargin{dimensionality information: how many countries, etc}

% \zhmargin{add more details on the setup}

The local response map
\(
f_i(x_{i,k},u_{i,k},x_{i,0})
\)
is product-wise defined by a Bass-type diffusion model~\cite{kiesling2012agent} with \(\vartheta\!\in\!\{\!1,\!2,\!3\}\):
\begin{equation*}
\begin{aligned}
f_i^{(\vartheta)}
&=
\langle w_i^{(\vartheta)}\!,z_{i,k}\rangle
+(1\!-\!z_{i,k}^{(\vartheta)})
\Big[
\psi_{i}^{(\vartheta)}\!(u_{i,k}^{(\vartheta)})
\!+\!\ell_{i}^{(\vartheta)}z_{i,k}^{(\vartheta)}
\Big],\\
z_{i,k}^{(\vartheta)}
&=
\lambda x_{i,k}^{(\vartheta)}
+(1\!-\!\lambda)x_{i,0}^{(\vartheta)},
\quad
x_{i,0}^{(\vartheta)}\sim\mu_{i,0}^{(\vartheta)}.
\end{aligned}
\end{equation*}
Here, the inner product \(\langle w_i^{(\vartheta)}\!,z_{i,k}\rangle\) represents the coupling effect among different products within agent \(i\), where \(w_{i}^{(\vartheta)}\in[-1,1]^3\). The coefficient \(\ell_{i}^{(\vartheta)} > 0\) represents the imitation
effect, and \(\lambda \in (0,1)\) balances the current and initial adoption
levels. 
% \zhmargin{check the parameter ranges}
The local states are coupled through the network and evolve according to an additive network structure
\begin{equation}
x_{k+1}
=
\alpha A x_k
+(1-\alpha)f(x_k,u_k,x_0),
\label{eq:ev_network_dynamics}
\end{equation}
where \(A \in \mathbb{R}^{3N \times 3N} \) is doubly stochastic and corresponds to a ring graph augmented with second-neighbor and long-range links, with randomly perturbed symmetric weights.
In \eqref{eq:ev_network_dynamics}, \(\alpha\in(0,1)\) controls cross-market coupling, and
\(
f=\operatorname{col}(f_1,\ldots,f_N).
\)
The objective for product \(\vartheta\!\in\!\{1,\!2,\!3\}\) at agent \(i\in\{1,\ldots,N\}\) is
\begingroup
\setlength{\abovedisplayskip}{1pt}
\setlength{\abovedisplayshortskip}{1pt}
\setlength{\belowdisplayskip}{1pt}
\setlength{\belowdisplayshortskip}{1pt}
\begin{equation}
\Phi_i^{(\vartheta)}
\!\!=
\!-\!\left(\!\nu_i^{(\vartheta)}\!\! u_i^{(\vartheta)}\!\!\!-c_i^{(\vartheta)}\!\right)\psi_i^{(\vartheta)}\!(\!u_i^{(\vartheta)}\!)
-\beta_i^{(\vartheta)}\!\!
\left(
\!\sqrt{\!x_{i,\mathrm{ss}}^{(\vartheta)}\!+\!\epsilon_{\mathrm{s}}}
\!-\!\sqrt{\epsilon_{\mathrm{s}}}
\right)\!,
\label{eq:results_objective}
\end{equation}
\endgroup
where \(\nu_i^{(\vartheta)}\!>\!0\) is a multiplicative factor that links the battery-price decision with the EV price, \(c_i^{(\vartheta)}\!>\!0\) is the cost, \(\beta_i^{(\vartheta)}\!>\!0\) is a weighting parameter for the adoption benefit, \(\epsilon_{\mathrm{s}}\!=\!10^{-6}\) regularizes the square-root term near zero, and $x_{i,\mathrm{ss}}$ is the steady-state adoption level of agent $i$ (i.e., $x_{\mathrm{ss}}$ is the fixed-point of \eqref{eq:ev_network_dynamics}). In \eqref{eq:results_objective}, the first term is the
negative local profit, while the second term rewards EV adoption. The
network-wide objective is
\(
\Phi(u,x_{\mathrm{ss}})
\!=\!
\sum_{i=1}^{N}\!\sum_{\vartheta}\!\Phi_i^{(\vartheta)}(u_i^{(\vartheta)},x_{i,\mathrm{ss}}^{(\vartheta)}).
\)
Hence, minimizing $\Phi(u,x_{\mathrm{ss}})$ seeks a balance between the profit and the long-run EV adoption level across the network. 
% The experimental model is specified by \eqref{eq:results_bass}, \eqref{eq:ev_network_dynamics}, and \eqref{eq:results_objective}.

% \subsection{Benchmark Algorithms}\label{sec:experiment_algs}
We compare five methods in total: C-NAD, the proposed
D-NAD, local NAD (L-NAD), performative gradient method \cite{perdomo2020performative}, and
zero-order gradient method \cite{duchi2015optimal}. 
C-NAD uses the complete network-response sensitivity. L-NAD is obtained by setting \(\mathcal S_i=\{i\}\), such that each agent constructs its response approximation without receiving derivatives from neighboring agents.
The performative gradient
\begingroup
\setlength{\abovedisplayskip}{2pt}
\setlength{\abovedisplayshortskip}{2pt}
\setlength{\belowdisplayskip}{2pt}
\setlength{\belowdisplayshortskip}{2pt}
\[
\widehat g_{i,k}^{\mathrm{perf}}
=
\frac{1}{N_{\mathrm{mb}}}
\sum_{m=1}^{N_{\mathrm{mb}}}
\partial_{u_i}\!\Phi_i(u_{i,k},x_{i,k}^m)
\]
\endgroup
retains the direct derivative and does not address the decision-induced response~\cite{perdomo2020performative}.
The zero-order gradient method using the standard two-point estimator~\cite{duchi2015optimal} is
\begingroup
\setlength{\abovedisplayskip}{2pt}
\setlength{\abovedisplayshortskip}{2pt}
\setlength{\belowdisplayskip}{2pt}
\setlength{\belowdisplayshortskip}{2pt}
\[
\widehat g_{i,k}^{\mathrm{ZO}}
\!=\!
\frac{n_uN}{2c_kN_{\mathrm{dir}}}
\sum_{s=1}^{N_{\mathrm{dir}}}
% \frac{d}{2c_k}
\left[
\widehat{\Phi}_k(u_k\!+\!c_k\delta_k^s)
\!-\!
\widehat{\Phi}_k(u_k\!-\!c_k\delta_k^s)
\right]
\!\delta_{i,k}^s,
\]
\endgroup
which uniformly samples an exploration direction $\delta_k^s$
% \(
% \delta_k^s
% \sim
% \operatorname{Unif}(\mathbb S)
% \)
from a unit sphere
\(
\mathbb S
\!:=\!
\left\{
\delta\!\in\!\mathbb R^{n_uN}\!:\!
\|\delta\|_2\!=\!1
\right\}.
\)
We use \(N_{\mathrm{dir}} \in \mathbb{N}_{+} \) exploration directions. 
Further, \(c_k\!>\!0\) is the perturbation radius, \(\delta_{i,k}^s\) is the component associated with agent \(i\), and
\(\widehat{\Phi}_k\) is the
mini-batch objective. 
% Both are compared with the proposed D-NAD method in the
% experiments.
% \zhmargin{mention how to select and tune hyperparameters}

% \subsection{Numerical Results}
We select the largest numerically stable step size from
a grid spanning \(0.01\) to \(0.10\) with increments of \(0.005\).
This gives fixed step sizes of \(\eta=0.025\) for the zero-order
method and \(\eta=0.10\) for the other four methods.
All the five methods use a mini-batch size of \(N_{\mathrm{mb}}=10\). For the zero-order gradient method, we use \(N_{\mathrm{dir}}=10\) exploration directions and a diminishing perturbation radius
\(c_k=\frac{0.04}{\sqrt{k+1}}\). 
For the Neumann order in the D-NAD method, we empirically choose \(R=3\) to balance approximation accuracy and local computational cost; see also Figure~\ref{fig:objective_comparison}(b) for more details.

For evaluations, let \(\Phi_k\) denote the transient objective evaluated
at the \(k\)-th iterate. 
% using the fixed evaluation samples. 
Let \(u^\star\) denote the optimal solution by solving the steady-state
problem in an offline and centralized manner, and define the optimal value
\(\Phi^\star:=\widetilde{\Phi}(u^\star)\).
Moreover, \(\Phi_T^{(R)}\) denotes the final transient objective obtained
by D-NAD with the Neumann order \(R\).

Figure~\ref{fig:objective_comparison} compares the transient objective
errors. Figure~\ref{fig:objective_comparison}(a) reports the optimality gaps
\(|\Phi_k-\Phi^\star|\) of all the methods with \(T=1000\), and
Figure~\ref{fig:objective_comparison}(b) illustrates
\(|\Phi_T^{(R)}-\Phi^\star|\) for \(R=1,\ldots,10\) and \(T=1000\).

\begingroup
\setlength{\abovedisplayskip}{1pt}
\setlength{\abovedisplayshortskip}{1pt}
\setlength{\belowdisplayskip}{1pt}
\setlength{\belowdisplayshortskip}{1pt}
\begin{figure}[!tb]
    \centering
    \setlength{\abovecaptionskip}{-3pt}
    \includegraphics[width=0.495\columnwidth]
    {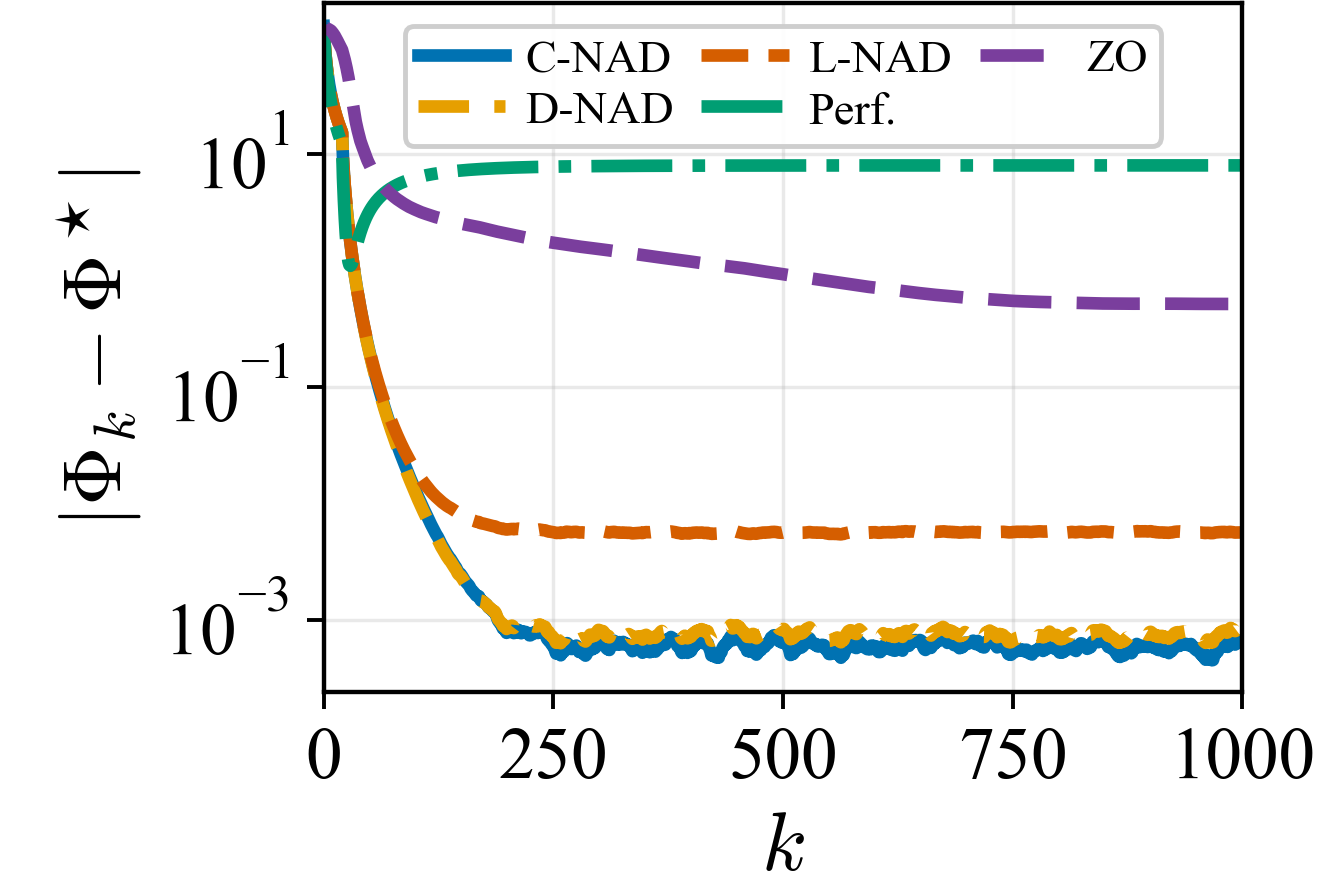}%
    \hspace{0.005\columnwidth}%
    \includegraphics[width=0.495\columnwidth]
    {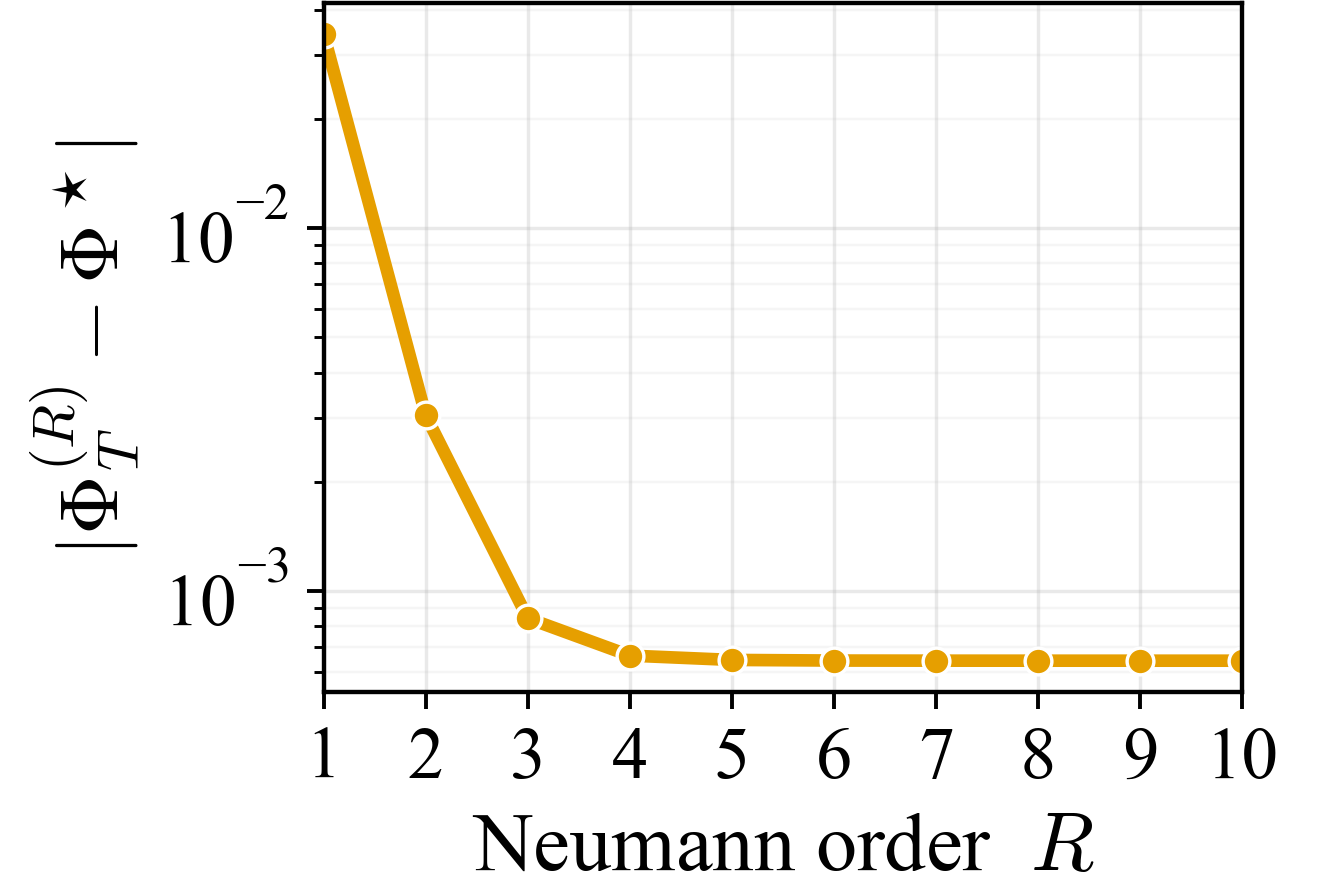}
    \caption{Objective comparison:
    (a) objective errors of different methods;
    (b) final transient D-NAD error versus the Neumann order $R$.}
    \label{fig:objective_comparison}
\end{figure}
\endgroup

Figure~\ref{fig:objective_comparison}(a) shows that D-NAD closely tracks C-NAD and achieves a smaller residual error than L-NAD, whereas the performative and ZO methods result in substantially larger errors.
Figure~\ref{fig:objective_comparison}(b) shows that the error decreases sharply from \(R=1\) to \(R=3\) and changes marginally for \(R\geq4\), supporting the use of \(R=3\) to balance accuracy and costs.

Figure~\ref{fig:decision_comparison} compares the price decisions for
plug-in hybrid electric vehicles. Figure~\ref{fig:decision_comparison}(a) demonstrates the distributions of agent-wise
final decisions, where the diamonds denote the means, and
$\sigma^2$ denotes the variance.
Figure~\ref{fig:decision_comparison}(b) shows the corresponding
agent-average decision trajectories. They indicate that C-NAD, D-NAD,
and L-NAD yield similar decisions, with only small differences.

\begingroup
\setlength{\abovedisplayskip}{1pt}
\setlength{\abovedisplayshortskip}{1pt}
\setlength{\belowdisplayskip}{1pt}
\setlength{\belowdisplayshortskip}{1pt}
\begin{figure}[!tb]
    \centering
    \setlength{\abovecaptionskip}{-3pt}
    \includegraphics[width=0.495\columnwidth]
    {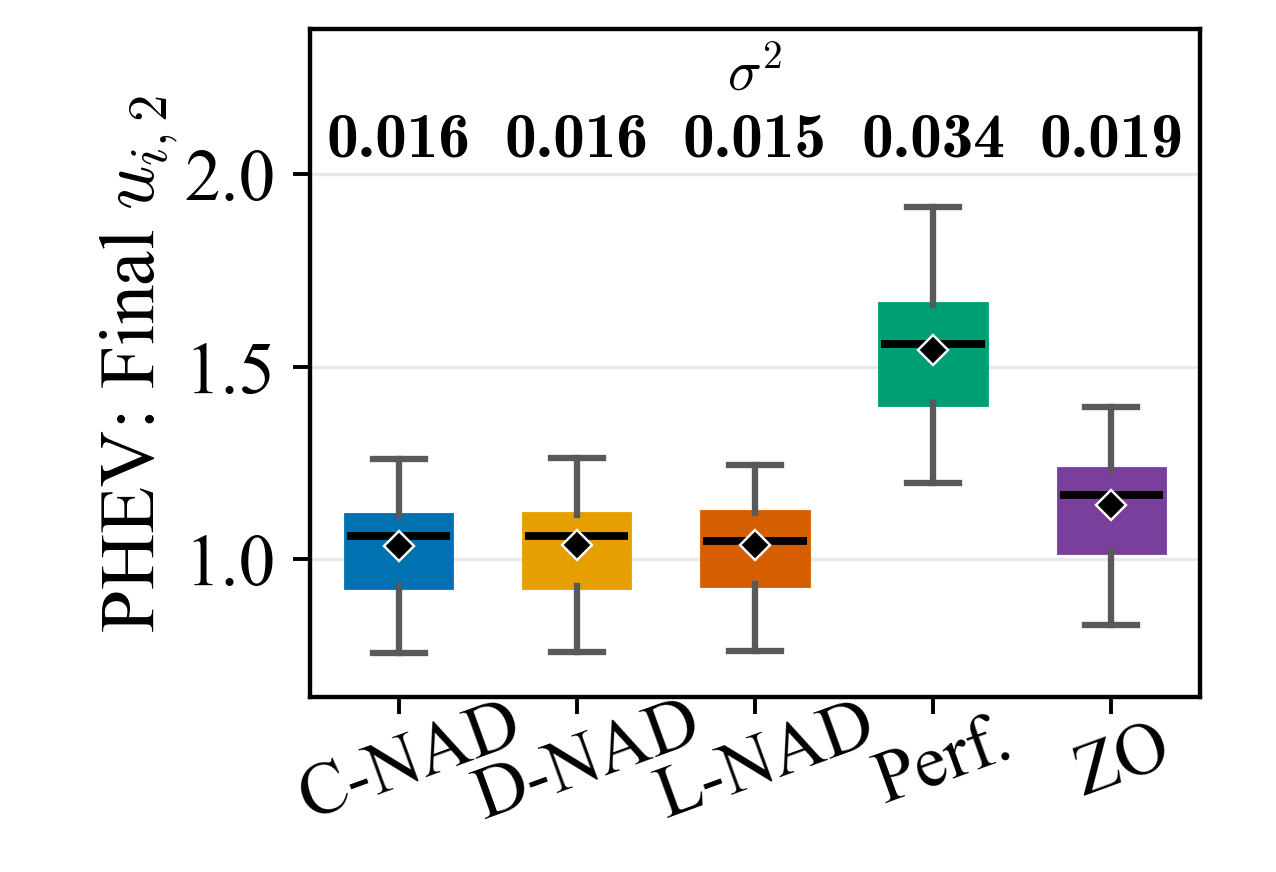}%
    \hspace{0.005\columnwidth}%
    \includegraphics[width=0.495\columnwidth]
    {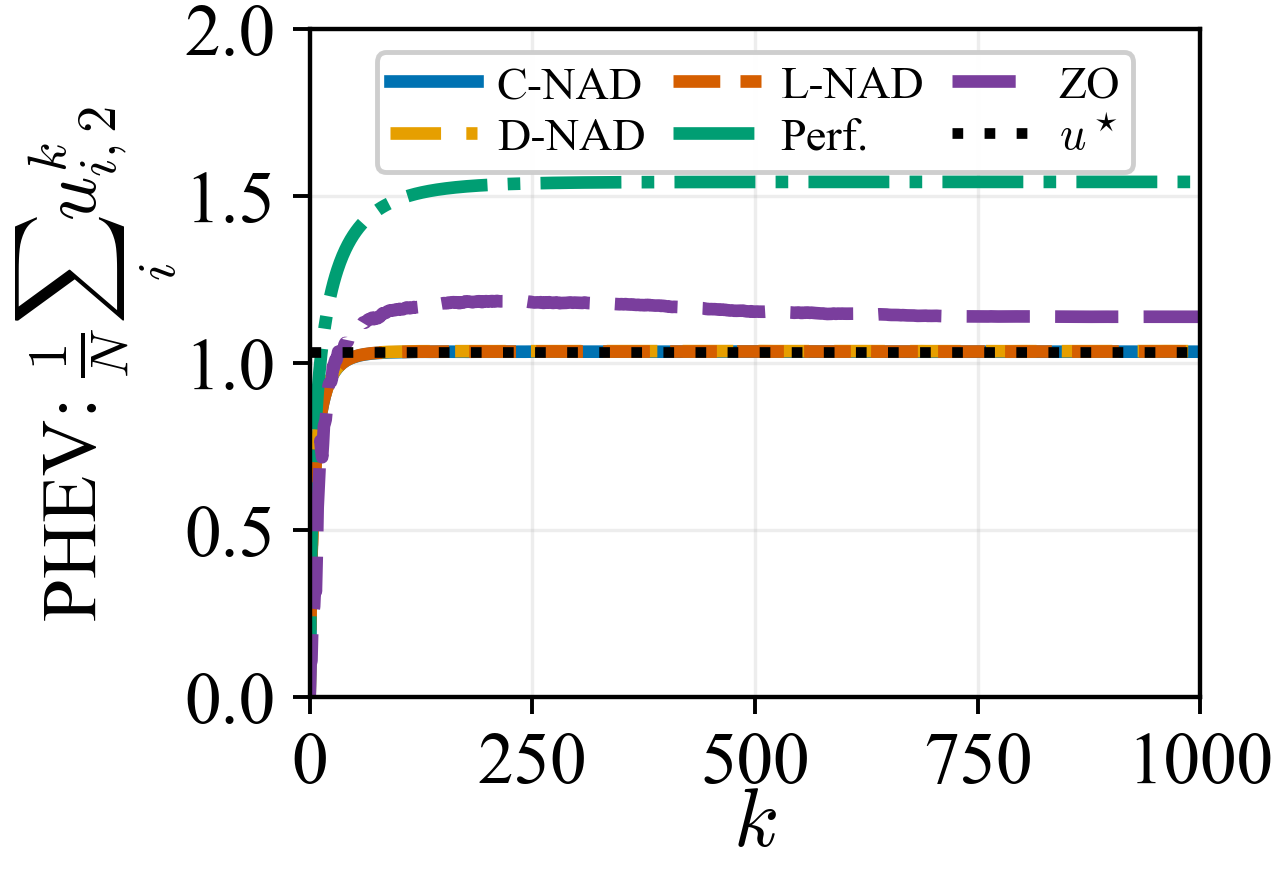}
    \caption{Decision comparison for plug-in hybrid electric vehicles:
    (a) final agent-wise decisions;
    (b) evolution of the agent-average decision.}
    \label{fig:decision_comparison}
\end{figure}
\endgroup

Figure~\ref{fig:decision_comparison} shows that the final decisions and
trajectory produced by D-NAD are close to those of C-NAD and the benchmark $u^\star$. The performative baseline converges to a
substantially larger decision because it ignores the decision-dependent state response, while ZO brings a larger residual error due to exploration noise.

\section{Conclusion}
\label{sec:conclusion}

We studied decision-dependent stochastic optimization with network distribution dynamics, where the collective decisions of agents induce local distribution shifts coupled through network interactions. We proposed the D-NAD algorithm, which approximates the network-response sensitivity using local Jacobian information, neighbor communication, and a truncated Neumann expansion. We characterized the closed-loop performance guarantee through the number of iterations and the approximation errors. We demonstrated the performance of our algorithm in an EV price--demand case study and compared it against competing benchmarks.

% We introduce the C-NAD gradient to capture
% both the direct and network-mediated effects of a decision, and develop
% D-NAD using local information and a truncated Neumann expansion. We
% characterize its approximation error and show through EV price--demand
% experiments that D-NAD remains close to C-NAD while outperforming
% L-NAD, performative and zero-order methods.
% This work studied decision-dependent stochastic optimization in
% networked dynamic systems. We introduced the network-aware
% decision-response (NAD) gradient, which accounts for both the direct
% effect of a decision and its indirect effect through the steady-state
% network response. We then proposed D-NAD, which approximates this
% response using local and neighboring information together with a
% truncated Neumann expansion.

% The D-NAD gradient error consists of a Neumann truncation error and a
% neighborhood approximation error, both of which depend on the Neumann
% order \(R\). For every finite \(R\), the total error is bounded under the
% stated assumptions. As \(R\to\infty\), the truncation contribution
% vanishes, while the neighborhood contribution remains uniformly
% bounded, yielding a bounded limiting error. The EV price--demand experiments
% showed that D-NAD remains close to centralized NAD and improves upon
% local-only and performative updates.
% \zhmargin{These two paragraphs can be simplified.}

Future directions include, but are not limited to 
% i) handling time-varying, directed network topology, 
i) developing constrained extensions of D-NAD, ii) addressing risk-sensitive objectives such as conditional value-at-risk to
hedge against tail outcomes, and iii) analyzing decision-dependent zero-sum network games with min-max formulations.

% work includes learning the response model and network structure
% online, particularly when the topology is uncertain or time varying.
% Another direction is on-demand sampling, where each agent adaptively
% % selects when and how many samples to acquire to balance estimation
% % accuracy and sampling cost. The framework can also be extended to
% % % decision-dependent zero-sum games with minimax objectives or to
% % % risk-sensitive optimization using conditional value-at-risk (CVaR) to
% % control adverse tail outcomes. Other directions include constrained and
% second-order extensions of D-NAD.

% \zhmargin{Check references; correct those arXiv preprints (some are published); add new references in 26.}

\balance
\bibliographystyle{IEEEtran}
\bibliography{bibliography}

\newpage
\nobalance % balancing is only meant for the bibliography page; left on, it distorts line spacing in the appendix
\appendix
% \section{Proofs}
% \label{app:proofs}

\subsection{Proof of Theorem~\ref{lem:average-neumann-error}}
\label{app:proof-neumann-error}

\begin{proof}
For each mini-batch sample \(m\), the tail of the
Neumann series can be written as
\(
\sum_{r=R+1}^{\infty}((H_k^m)^\top)^r
\!=\!
((H_k^m)^\top)^{R+1}(I\!-\!H_k^m)^{-\top}.
\)
Therefore,
\[
e_{k}^{(1)}\!(R)
=
\tfrac{1}{N_{\mathrm{mb}}}
{\textstyle\sum\nolimits_{m=1}^{N_{\mathrm{mb}}}
(B_k^m)^\top((H_k^m)^\top)^{R+1}
(I-H_k^m)^{-\top}q_k^m}.
\]
The triangle inequality and the sub-multiplicative property of the matrix norm give
\[
\big\|
e_{k}^{(1)}\!(R)
\big\|
\!\leq\!
\tfrac{1}{N_{\mathrm{mb}}}
\!{\textstyle\sum\nolimits_{m=1}^{N_{\mathrm{mb}}}
\!\|B_k^m\|\!
\|H_k^m\|^{R+1}\!
\|(I\!-\!H_k^m)^{-1}\|\!
\|q_k^m\|}.
\]
By Assumption~\ref{asm:network-dynamics} and \eqref{eq:gamma_k_definition}, we have \(\|H_k^m\|\!\leq\!\gamma_k\!<\!1\) and \(\|B_k^m\|\leq(1-\alpha)L_f^u\). With
\(
\|(I-H_k^m)^{-1}\|
\leq
\sum_{r=0}^{\infty}\|H_k^m\|^r=
\frac{1}{1-\|H_k^m\|}
\leq
\frac{1}{1-\gamma_k},
\)
we obtain
\[
\big\|
e_{k}^{(1)}(R)
\big\|^2
\leq
C_0
\tfrac{\gamma_k^{2R+2}}
{(1-\gamma_k)^2},
\tag{A.1}\label{eq:app-neumann-error}
\]
where Assumption~\ref{asm:objective-regularity} gives
\(
\|q_k^m\|\leq L_\Phi^x
\)
for every \(m\).
Taking expectations in \eqref{eq:app-neumann-error} and averaging over \(k\) proves the result.
\end{proof}

\subsection{Proof of Theorem~\ref{lem:average-neighborhood-error}}
\label{app:proof-neighborhood-error}

\begin{proof}
By definitions of \(\Delta_i\) and \eqref{eq:extended_H}, for every sample \(m\),
\(
H_k^m\!-\!\widetilde H_{\mathcal S_i,k}^m
\!=\!
\alpha\Delta_i.
\)
Using the telescoping identity
\(
X^r\!-\!Y^r
\!=\!
\sum_{\ell=0}^{r-1}\!
X^{r-1-\ell}(X\!-\!Y)Y^\ell,
\)
we have
\begin{align*}
\|\!(\!H_k^m)^r\!\!-\!(\!\widetilde H_{\mathcal S_i,k}^m\!)^r\!\|
&\!\leq\!
{\textstyle\sum\nolimits_{\ell=0}^{r-1}}
\|\!H_k^m\!\|^{r\!-\!1\!-\ell}
\|\!H_k^m\!-\!\widetilde H_{\mathcal S_i,k}^m\!\|
\|\!\widetilde H_{\mathcal S_i,k}^m\!\|^\ell\\
&\!\leq\!
\alpha r\|\Delta_i\|\gamma_k^{r-1}.
\tag{A.2}\label{eq:app-matrix-power}
\end{align*}
Applying the triangle inequality and \eqref{eq:app-matrix-power} gives
\[
\|e_{i,k}^{(2)}(R)\|
\!\leq\!
\tfrac{\alpha\|\Delta_i\|}{N_{\mathrm{mb}}}
{\textstyle\sum\nolimits_{m=1}^{N_{\mathrm{mb}}}
\|B_{ii,k}^m\|\|q_k^m\|
\cdot\sum_{r=2}^{R}\!r\gamma_k^{r-1}}.
\]
Assumptions 1 and 2 give
\(
\|B_{ii,k}^m\|\!\leq\!(1\!-\!\alpha)L_f^u
\)
and
\(
\|q_k^m\|\!\leq\! L_\Phi^x
\)
for every \(m\).
Thus, we obtain
\(
\big\|
e_{i,k}^{(2)}(R)
\big\|^2
\!\leq\!
\alpha^2C_0
\|\Delta_i\|^2
\kappa_R(\gamma_k)^2.
\)
Summing over all agents gives
\(
\big\|
e_{k}^{(2)}\!(R)\!
\big\|^2
\!\leq\!
\alpha^2C_0
\Delta_A
\kappa_R(\gamma_k)^2.
\)
Taking expectations and averaging the bound on \(\|e_k^{(2)}(R)\|^2\) derived from \eqref{eq:app-matrix-power} give
\(
\frac1T\sum_{k=0}^{T-1}\mathbb E[\|e_k^{(2)}\|^2]
\leq
C_0\mathcal B^{(2)}.
\)
Finally, since \(\gamma_k<1\),
\[
\kappa_R(\gamma_k)
\leq
{\textstyle\sum\nolimits_{r=2}^{\infty}}r\gamma_k^{r-1}
=
\tfrac{1}{(1-\gamma_k)^2}-1
=
\kappa_\infty(\gamma_k).
\tag{A.3}\label{eq:app-kappa-bound}
\]
Applying \eqref{eq:app-kappa-bound} to the definition of \(\mathcal B^{(2)}\) proves
\(
\mathcal B^{(2)}
\leq
\mathcal B^{(2)}_\infty.
\)
\end{proof}

\subsection{Proof of Theorem~\ref{thm:distributed-nad-convergence}}
\label{app:proof-dnad-convergence}

To analyze convergence, we define the gradient errors \(e_k^{\mathrm{C\text{-}NAD}}
\!\!:=\!
\nabla\widetilde{\Phi}(u_k)\!-\!\widehat g_k^{\mathrm{C\text{-}NAD}}
\)
and
\(e_k^{\mathrm{D\text{-}NAD}}
\!\!:=\!
\nabla\widetilde{\Phi}(u_k)\!-\!\widehat g_k^{\mathrm{D\text{-}NAD}}\!.
\) 
% through
% \(
% \nabla\widetilde{\Phi}(u_k)
% =
% \widehat g_k^{\mathrm{C\text{-}NAD}}
% +
% e_k^{\mathrm{C\text{-}NAD}}
% =
% \widehat g_k^{\mathrm{D\text{-}NAD}}
% +
% e_k^{\mathrm{D\text{-}NAD}}.
% \)
For a fixed \(R\), let
\(
e_k^{\mathrm{appr}}:=e_k^{(1)}(R)+e_k^{(2)}(R)
\)
denote the combined approximation
error.
The D-NAD gradient error satisfies
\(
e^{\mathrm{D\text{-}NAD}}_{k}
\!=\!
e_k^{\mathrm{appr}}\!+\!e^{\mathrm{C\text{-}NAD}}_{k}.
\)
Write
\(\mathbb E_k[\cdot]:=\mathbb E[\cdot\mid\!\mathcal F_{k-1}]\) and
\(\operatorname{Var}_k[\cdot]:=\operatorname{Var}[\cdot\mid\!\mathcal F_{k-1}]\) for each \(k\). Before proving Theorem~\ref{thm:distributed-nad-convergence}, we present the following lemma that bounds the conditional variance of the D-NAD gradient estimate.

\begin{lemma}[D-NAD sampling variance]
\label{cor:dnad-variance}
Suppose Assumption~\ref{asm:gradient-variance} holds and the mini-batch
samples are conditionally independent given \(\mathcal F_{k-1}\).
Then
\begin{align*}
\operatorname{Var}_k\!\left[
\widehat g_k^{\mathrm{D\text{-}NAD}}
\right]
&\le
\tfrac{4M_V}{N_{\mathrm{mb}}}
\left\|\mathbb E_k\!\left[
\widehat g_k^{\mathrm{D\text{-}NAD}}
\right]\right\|^2+\tfrac{2M}{N_{\mathrm{mb}}}\\
&\qquad+
\left(2+\tfrac{4M_V}{N_{\mathrm{mb}}}\right)
\mathbb E_k\!\left[
\|e_k^{\mathrm{appr}}\|^2
\right].
\tag{A.4}\label{eq:app-dnad-variance}
\end{align*}
\end{lemma}

\begin{proof}
Since
\(\widehat g_k^{\mathrm{C\text{-}NAD}}
=\widehat g_k^{\mathrm{D\text{-}NAD}}+e_k^{\mathrm{appr}}\),
taking conditional expectations gives
\begin{align*}
\|\bar g_k^{\mathrm{C\text{-}NAD}}\|^2
&=
\Bigl\|
\mathbb E_k\!\left[\widehat g_k^{\mathrm{D\text{-}NAD}}\right]
+\mathbb E_k\!\left[e_k^{\mathrm{appr}}\right]
\Bigr\|^2\\
&\overset{\mathrm{(a.1)}}{\le}
2\Bigl\|\mathbb E_k\!\left[\widehat g_k^{\mathrm{D\text{-}NAD}}\right]\Bigr\|^2
+2\Bigl\|\mathbb E_k\!\left[e_k^{\mathrm{appr}}\right]\Bigr\|^2\\
&\overset{\mathrm{(a.2)}}{\le}
2\Bigl\|\mathbb E_k\!\left[\widehat g_k^{\mathrm{D\text{-}NAD}}\right]\Bigr\|^2
+2\,\mathbb E_k\!\left[\|e_k^{\mathrm{appr}}\|^2\right],
\end{align*}
where (a.1) uses \(\|a+b\|^2\le 2\|a\|^2+2\|b\|^2\) and
(a.2) follows from conditional Jensen's inequality in the form \(\|\mathbb E_k[\cdot]\|^2\!\leq\!\mathbb E_k[\|\!\cdot\!\|^2]\).
The conditional variance satisfies
\begin{align*}
\operatorname{Var}_k\!\left[
\widehat g_k^{\mathrm{D\text{-}NAD}}\right]
&=
\operatorname{Var}_k\!\left[
\widehat g_k^{\mathrm{C\text{-}NAD}}-e_k^{\mathrm{appr}}\right]\\
&\overset{\mathrm{(a.3)}}{\le}
2\operatorname{Var}_k\!\left[
\widehat g_k^{\mathrm{C\text{-}NAD}}\right]
+2\operatorname{Var}_k\!\left[e_k^{\mathrm{appr}}\right]\\
&\overset{\mathrm{(a.4)}}{\le}
\tfrac{2(M+M_V\|\bar g_k^{\mathrm{C\text{-}NAD}}\|^2)}
{N_{\mathrm{mb}}}
+2\,\mathbb E_k\!\left[\|e_k^{\mathrm{appr}}\|^2\right]\\
&\overset{\mathrm{(a.5)}}{\le}
\tfrac{4M_V}{N_{\mathrm{mb}}}
\Bigl\|\mathbb E_k\!\left[\widehat g_k^{\mathrm{D\text{-}NAD}}\right]\Bigr\|^2
+\tfrac{2M}{N_{\mathrm{mb}}}\\
&\qquad+
\left(2+\tfrac{4M_V}{N_{\mathrm{mb}}}\right)
\mathbb E_k\!\left[\|e_k^{\mathrm{appr}}\|^2\right].
\end{align*}
In (a.3), we apply the inequality in (a.1) to the conditionally centered
terms. For (a.4), with Assumption~\ref{asm:gradient-variance}, conditional independence makes the variance of the
mini-batch average equal to the sum of its sample variances divided
by \(N_{\mathrm{mb}}^2\). We bound the conditional
variance of \(e_k^{\mathrm{appr}}\) by its conditional second moment.
Finally, (a.5) substitutes the bound in (a.2).
\end{proof}

\begin{proof}[Proof of Theorem~\ref{thm:distributed-nad-convergence}]
% Keep the explicit label: this proof follows the auxiliary corollary and its proof.
By \(L_{\widetilde \Phi}\)-smoothness and the D-NAD update \(u_{k+1}\!-\!u_k=-\eta\,\widehat g_k^{\mathrm{D\text{-}NAD}}\), we have
\begin{align*}
\widetilde\Phi(u_{k+1})
&\le
\widetilde\Phi(u_k)-\eta\bigl\langle\nabla\widetilde\Phi(u_k),
\widehat g_k^{\mathrm{D\text{-}NAD}}\bigr\rangle\\
&\quad+
\tfrac{L_{\widetilde\Phi}\eta^2}{2}
\|\widehat g_k^{\mathrm{D\text{-}NAD}}\|^2.
\tag{A.5}\label{eq:app-smooth-descent}
\end{align*}
By the conditional variance decomposition and
\eqref{eq:app-dnad-variance},
\begin{align*}
\mathbb E_k\!\left[
\|\widehat g_k^{\mathrm{D\text{-}NAD}}\|^2
\right]
&=
\operatorname{Var}_k\!\left[\widehat g_k^{\mathrm{D\text{-}NAD}}\right]
+\left\|\mathbb E_k\!\left[
\widehat g_k^{\mathrm{D\text{-}NAD}}
\right]\right\|^2\\
&\le
\left(1+\tfrac{4M_V}{N_{\mathrm{mb}}}\right)
\left\|\mathbb E_k\!\left[
\widehat g_k^{\mathrm{D\text{-}NAD}}
\right]\right\|^2
+\tfrac{2M}{N_{\mathrm{mb}}}\\
&\quad+
\left(2+\tfrac{4M_V}{N_{\mathrm{mb}}}\right)
\mathbb E_k\!\left[\|e_k^{\mathrm{appr}}\|^2\right].
\tag{A.6}\label{eq:app-conditional-moment}
\end{align*}
Substituting \eqref{eq:app-conditional-moment} into the conditional expectation of \eqref{eq:app-smooth-descent} and using \(
\mathbb E_k\!\left[\widehat g_k^{\mathrm{D\text{-}NAD}}\right]\!=\!\nabla\widetilde\Phi(u_k)
\!-\!\mathbb E_k\!\left[e_k^{\mathrm{D\text{-}NAD}}\right],
\) gives
\begin{align*}
\mathbb E_k\!\left[\widetilde\Phi(u_{k+1})\right]
&\le
\widetilde\Phi(u_k)-\eta\bigl\langle
\nabla\widetilde\Phi(u_k),
\mathbb E_k\!\left[\widehat g_k^{\mathrm{D\text{-}NAD}}\right]
\bigr\rangle\\
&\;+
\tfrac{L_{\widetilde\Phi}\eta^2}{2}
\left(1+\tfrac{4M_V}{N_{\mathrm{mb}}}\right)
\left\|\mathbb E_k\!\left[
\widehat g_k^{\mathrm{D\text{-}NAD}}
\right]\right\|^2\\
&\;+
L_{\widetilde\Phi}\eta^2\!
\left(1+\tfrac{2M_V}{N_{\mathrm{mb}}}\right)\!\mathbb E_k\!\left[\|e_k^{\mathrm{appr}}\|^2\right]+\tfrac{L_{\widetilde\Phi}M\eta^2}{N_{\mathrm{mb}}}\\
&\overset{\mathrm{(a.6)}}{\le}\widetilde\Phi(u_k)
-C_2\|\nabla\widetilde\Phi(u_k)\|^2
+C_4\\
&\hspace{-3em}+C_1\bigl\|
\mathbb E_k[e_k^{\mathrm{D\text{-}NAD}}]
\bigr\|^2+
C_3\,\mathbb E_k\!\left[\|e_k^{\mathrm{appr}}\|^2\right],
\tag{A.7}\label{eq:app-conditional-descent}
\end{align*}
where (a.6) uses \(2\langle a,b\rangle\!\le\!\|a\|^2\!+\!\|b\|^2\)
and \(\|a\!-\!b\|^2\!\le\!2\|a\|^2\!+\!2\|b\|^2\), with
\begin{alignat*}{2}
C_1&:=\tfrac{\eta}{2}+L_{\widetilde\Phi}\eta^2
\bigl(1+\tfrac{4M_V}{N_{\mathrm{mb}}}\bigr),
&\quad
C_3&:=L_{\widetilde\Phi}\eta^2
\bigl(1+\tfrac{2M_V}{N_{\mathrm{mb}}}\bigr),\\
C_2&:=\tfrac{\eta}{2}-L_{\widetilde\Phi}\eta^2
\bigl(1+\tfrac{4M_V}{N_{\mathrm{mb}}}\bigr),
&\quad
C_4&:=\tfrac{L_{\widetilde\Phi}M\eta^2}{N_{\mathrm{mb}}}.
\tag{A.8}\label{eq:app-descent-constants}
\end{alignat*}
Taking conditional expectations in the error decomposition
\(e_k^{\mathrm{D\text{-}NAD}}
\!=\!e_k^{\mathrm{C\text{-}NAD}}\!+\!e_k^{\mathrm{appr}}\)
and applying \(\|a\!+\!b\|^2\!\le\!2\|a\|^2\!+\!2\|b\|^2\) and Jensen's inequality give
\[
\|\mathbb E_k[e_k^{\mathrm{D\text{-}NAD}}]\|^2
\leq
2\|\mathbb E_k[e_k^{\mathrm{C\text{-}NAD}}]\|^2
+2\,\mathbb E_k[\|e_k^{\mathrm{appr}}\|^2].
\tag{A.9}\label{eq:app-error-split}
\]
Substituting \eqref{eq:app-error-split} into \eqref{eq:app-conditional-descent},
taking expectations, averaging over \(k\), and using \(\widetilde\Phi^\star:=\inf_u\widetilde\Phi(u)
\leq\widetilde\Phi(u_T)\) yield
\begin{align*}
\tfrac{C_2}{T}{\textstyle\sum\nolimits_{k=0}^{T-1}}
&\mathbb E\bigl[\|\nabla\widetilde\Phi(u_k)\|^2\bigr]
\le
\tfrac{\widetilde\Phi(u_0)-\widetilde\Phi^\star}{T}
+C_4\\
&\quad+\tfrac{2C_1+C_3}{T}{\textstyle\sum\nolimits_{k=0}^{T-1}}\mathbb E\bigl[\|e_k^{\mathrm{appr}}\|^2\bigr]\\
&\quad+
\tfrac{2C_1}{T}{\textstyle\sum\nolimits_{k=0}^{T-1}}\mathbb E\bigl[
\|\mathbb E_k[e_k^{\mathrm{C\text{-}NAD}}]\|^2
\bigr].
\tag{A.10}\label{eq:app-averaged-descent}
\end{align*}

To bound the C-NAD term in \eqref{eq:app-averaged-descent},
fix the realized decision sequence and define
\[
V_k:=
\mathbb E_{(x_0,d)}\!\left[\|x_k-h_A(u_k,d)\|\right],
\tag{A.11}\label{eq:app-tracking-definition}
\]
where the expectation is over the joint distribution of \((x_0,d)\).
The transient state \(x_k\) and
\(x_{\mathrm{ss}}=h_A(u_k,d)\), which is its steady-state counterpart, use the same \(d\).

Differentiating the fixed-point equation for \(h_A\) shows that
the summand in \eqref{eq:centralized_gd}, evaluated at
\(x=h_A(u_k,d)\), equals
\(\left.\nabla_u[\Phi(u,h_A(u,d))]\right|_{u=u_k}\).
The mean of the steady-state summand is \(\nabla\widetilde\Phi(u_k)\) by
Assumption~\ref{asm:lipschitz-domination}. Provided that
\(\mathbb E_k[\widehat g_k^{\mathrm{C\text{-}NAD}}]\) equals the
mean of the summand evaluated at the transient state
\(x_k\), Jensen's inequality and Lipschitz continuity in \(x\) give
\[
\|\mathbb E_k[e_k^{\mathrm{C\text{-}NAD}}]\|=\|\nabla \widetilde\Phi(u_k)\!-\!\mathbb E_k[\widehat g_k^{\mathrm{C\text{-}NAD}}]\|
\leq L_C^x V_k,
\tag{A.12}\label{eq:app-cnad-bias}
\]
where
\(
L_C^x=
M_\Phi^u
+\frac{L_f^uM_\Phi^x+L_\Phi^xM_f^u}{1-L_f^x}
+\frac{L_\Phi^xL_f^uM_f^x}{(1-L_f^x)^2}.
\)

We adapt the tracking argument in 
\cite[Lemma~19]{he2025decisiondependent} to the D-NAD iterates. By
contraction of \(f_A(\cdot,u_k,d)\) and the fixed-point identity
for \(h_A(u_k,d)\),
\begin{align*}
\|x_{k+1}-h_A(u_{k+1},d)\|
&\leq
\bigl[\alpha+(1\!-\!\alpha)L_f^x\bigr]
\|x_k\!-\!h_A(u_k,d)\|\\
&\hspace{-4em}+\|h_A(u_k,d)-h_A(u_{k+1},d)\|.
\tag{A.13}\label{eq:app-state-contraction}
\end{align*}
Subtracting the fixed-point equations at two decisions shows
that \(h_A\) is \(\frac{L_f^u}{1-L_f^x}\)-Lipschitz in \(u\).
Averaging \eqref{eq:app-state-contraction} over \((x_0,d)\) using
\(V_k\) from \eqref{eq:app-tracking-definition} and substituting
\(u_{k+1}\!-\!u_k\!=\!-\eta\,\widehat g_k^{\mathrm{D\text{-}NAD}}\)
give
\[
V_{k+1}
\leq
\bigl[\alpha+(1-\alpha)L_f^x\bigr]V_k
+\frac{\eta L_f^u}{1-L_f^x}
\|\widehat g_k^{\mathrm{D\text{-}NAD}}\|.
\tag{A.14}\label{eq:app-tracking-recursion}
\]
We then apply the weighted Cauchy--Schwarz inequality
\((\rho a+b)^2\leq\rho a^2+\frac{b^2}{(1-\rho)}\),
with \(\rho=\alpha+(1-\alpha)L_f^x\in(0,1)\) to \eqref{eq:app-tracking-recursion}, and obtain
\begin{align*}
V_{k+1}^2
&\leq
\bigl[\alpha+(1-\alpha)L_f^x\bigr]V_k^2\\
&\quad+\tfrac{\eta^2(L_f^u)^2}{(1-\alpha)(1-L_f^x)^3}
\|\widehat g_k^{\mathrm{D\text{-}NAD}}\|^2.
\tag{A.15}\label{eq:app-tracking-square}
\end{align*}
Subtracting \(V_k^2\) from both sides of \eqref{eq:app-tracking-square},
taking total expectations, and summing over \(k\) yield
\begin{align*}
&(1-\alpha)(1-L_f^x){\textstyle\sum\nolimits_{k=0}^{T-1}\mathbb E[V_k^2]}\\
&\leq
\!V_0^2\!-\!\mathbb E[V_T^2]\!+\!\tfrac{\eta^2(L_f^u)^2}{(1-\alpha)(1-L_f^x)^3}
\!{\textstyle\sum\nolimits_{k=0}^{T-1}
\!\mathbb E[\|\widehat g_k^{\mathrm{D\text{-}NAD}}\|^2]}.
\tag{A.16}\label{eq:app-tracking-telescoping}
\end{align*}
We drop the terminal contribution $\mathbb E[V_T^2]$ in \eqref{eq:app-tracking-telescoping},
divide by \((1\!-\!\alpha)(1\!-\!L_f^x)\), and apply the square of
\eqref{eq:app-cnad-bias}. This gives
\begin{align*}
&{\textstyle\sum\nolimits_{k=0}^{T-1}\mathbb E\!\left[
\|\mathbb E_k[e_k^{\mathrm{C\text{-}NAD}}]\|^2
\right]}\\
&\leq
\!\tfrac{(L_C^x)^2V_0^2}{(1\!-\!\alpha)(1\!-\!L_f^x)}
\!+\!\tfrac{\eta^2(L_C^xL_f^u)^2}{(1\!-\!\alpha)^2(1\!-\!L_f^x)^4}
\!{\textstyle\sum\nolimits_{k=0}^{T-1}
\!\mathbb E[\|\widehat g_k^{\mathrm{D\text{-}NAD}}\|^2]}.
\tag{A.17}\label{eq:app-bias-tracking}
\end{align*}

To bound the D-NAD second-moment sum in \eqref{eq:app-bias-tracking}, write
\[
\mathbb E_k[\widehat g_k^{\mathrm{D\text{-}NAD}}]
=\nabla\widetilde\Phi(u_k)
-\mathbb E_k[e_k^{\mathrm{C\text{-}NAD}}]
-\mathbb E_k[e_k^{\mathrm{appr}}].
\tag{A.18}\label{eq:app-mean-split}
\]
Substituting \eqref{eq:app-mean-split} into
\eqref{eq:app-conditional-moment} and applying
\(\|a+b+c\|^2\leq3(\|a\|^2+\|b\|^2+\|c\|^2)\)
together with conditional Jensen's inequality
\(\|\mathbb E_k[\cdot]\|^2
\leq\mathbb E_k[\|\cdot\|^2]\) give
\begin{align*}
&\mathbb E_k[\|\widehat g_k^{\mathrm{D\text{-}NAD}}\|^2]\\
&\quad\leq
3\left(1+\tfrac{4M_V}{N_{\mathrm{mb}}}\right)
\left(
\|\nabla\widetilde\Phi(u_k)\|^2
+\|\mathbb E_k[e_k^{\mathrm{C\text{-}NAD}}]\|^2
\right)\\
&\qquad+
\left(5+\tfrac{16M_V}{N_{\mathrm{mb}}}\right)
\mathbb E_k[\|e_k^{\mathrm{appr}}\|^2]
+\tfrac{2M}{N_{\mathrm{mb}}}.
\tag{A.19}\label{eq:app-conditional-moment-expanded}
\end{align*}
We then take expectations in \eqref{eq:app-conditional-moment-expanded}, sum over \(k\), and use
\(\mathbb E[\|e_k^{\mathrm{appr}}\|^2]
\!\leq\!2\!\left(\mathbb E[\|e_k^{(1)}\|^2]\!+\!\mathbb E[\|e_k^{(2)}\!\|^2]\right)\) to obtain
\begin{align*}
&{\textstyle\sum\nolimits_{k=0}^{T-1}
\mathbb E[\|\widehat g_k^{\mathrm{D\text{-}NAD}}\|^2]}
\!\leq\!
3\left(1+\tfrac{4M_V}{N_{\mathrm{mb}}}\right)
\!{\textstyle\sum\nolimits_{k=0}^{T-1}\mathbb E[
\|\nabla\widetilde\Phi(u_k)\|^2]}\\
&\quad+
3\left(1+\tfrac{4M_V}{N_{\mathrm{mb}}}\right)
{\textstyle\sum\nolimits_{k=0}^{T-1}\mathbb E\bigl[
\|\mathbb E_k[e_k^{\mathrm{C\text{-}NAD}}]\|^2
\bigr]}
+\tfrac{2MT}{N_{\mathrm{mb}}}\\
&\quad+
\left(10+\tfrac{32M_V}{N_{\mathrm{mb}}}\right)
{\textstyle\sum\nolimits_{k=0}^{T-1}\mathbb E\bigl[\|e_k^{(1)}\|^2+\|e_k^{(2)}\|^2\bigr]}.
\tag{A.20}\label{eq:app-moment-sum}
\end{align*}
Choose the step size \(\eta=\frac{\eta_0}{\sqrt T}\), where
\[
0\!<\!\eta_0\!\leq\!
\min\!\left\{
\tfrac{1}{4L_{\widetilde\Phi}
\left(1\!+\!\frac{4M_V}{N_{\mathrm{mb}}}\right)},
\tfrac{(1-\alpha)(1-L_f^x)^2}
{L_C^xL_f^u\sqrt{39\left(1\!+\!\frac{4M_V}{N_{\mathrm{mb}}}\right)}}
\right\}.
\tag{A.21}\label{eq:app-stepsize}
\]
The first restriction in \eqref{eq:app-stepsize} and the definitions
in \eqref{eq:app-descent-constants} give \(C_2\!\geq\!\frac{\eta}{4}\)
and \(C_1\!\leq\!\frac{3\eta}{4}\). The second restriction in
\eqref{eq:app-stepsize} gives
\(
\tfrac{\eta^2(L_C^xL_f^u)^2}{(1-\alpha)^2(1-L_f^x)^4}
\leq
\tfrac{1}{39\left(1+\frac{4M_V}{N_{\mathrm{mb}}}\right)T},
% \tag{A.22}\label{eq:app-tracking-multiplier}
\)
with which substituting \eqref{eq:app-moment-sum} into \eqref{eq:app-bias-tracking} gives a coefficient of at most \(\frac{1}{13T}\leq\frac{1}{13}\) for each of the
squared-gradient and squared conditional-mean sums. Moving the
squared conditional-mean sum to the left and dividing by
\(\frac{12}{13}\), we obtain
\begin{align*}
&{\textstyle\sum\nolimits_{k=0}^{T-1}
\!\mathbb E[\|\mathbb E_k[e_k^{\mathrm{C\text{-}NAD}}]\|^2]}
\!\!\leq\!\!\tfrac{1}{12}{\textstyle\sum\nolimits_{k=0}^{T-1}
\!\mathbb E[\|\nabla\widetilde\Phi(u_k)\|^2]}\\
&\qquad\quad+\!C_5\!+\!\tfrac{5}{18T}{\textstyle\sum\nolimits_{k=0}^{T-1}\!\!\left(\mathbb E[\|e_k^{(1)}\|^2]\!+\!\mathbb E[\|e_k^{(2)}\!\|^2]\right)},
\tag{A.22}\label{eq:app-bias-sum}
\end{align*}
where
\(
C_5=
\frac{13(L_C^x)^2V_0^2}{12(1-\alpha)(1-L_f^x)}
+\frac{M}{18N_{\mathrm{mb}}}.
\)

In \eqref{eq:app-bias-sum}, we use
\(10+\frac{32M_V}{N_{\mathrm{mb}}}
\leq10\left(1+\frac{4M_V}{N_{\mathrm{mb}}}\right)\)
and \(1+\frac{4M_V}{N_{\mathrm{mb}}}\geq1\) to bound the
approximation and sampling contributions, respectively.
Substituting \eqref{eq:app-bias-sum} into \eqref{eq:app-averaged-descent}
gives a squared-gradient coefficient of
\(C_2\!-\!\frac{C_1}{6}\!\geq\!\frac{\eta}{8}\) on the left.
Using
\[\tfrac{1}{T}{\textstyle\sum\nolimits_{k=0}^{T-1}\!\!\left(\mathbb E[\|e_k^{(1)}\|^2]\!+\!\mathbb E[\|e_k^{(2)}\!\|^2]\right)}
\!\leq\! C_0(\mathcal B^{(1)}+\mathcal B^{(2)})\] from
Theorems~\ref{lem:average-neumann-error}
and~\ref{lem:average-neighborhood-error}, and dividing by \(\frac{\eta}{8}\) yield
\begin{align*}
&\frac{1}{T}\!\sum_{k=0}^{T-1}
\mathbb E[\|\nabla\widetilde\Phi(u_k)\|^2]
\leq
\underbrace{
\frac{8(\widetilde\Phi(u_0)\!-\!\widetilde\Phi^\star)}{\eta T}
}_{\mathcal O(\frac{1}{\sqrt T})}
+
\underbrace{\frac{8C_4}{\eta}}
_{\mathcal O(\frac{1}{\sqrt T})}
+
\underbrace{\frac{16C_1C_5}{\eta T}}
_{\mathcal O(\frac{1}{T})}\\
&\quad+
\underbrace{
\frac{8C_0}{\eta}
\left(4C_1+2C_3+\frac{5C_1}{9T}\right)
\bigl(\mathcal B^{(1)}+\mathcal B^{(2)}\bigr)
}_{\mathcal O(\mathcal B^{(1)}+\mathcal B^{(2)})}.
\tag{A.23}\label{eq:app-normalized-descent}
\end{align*}
The orders in \eqref{eq:app-normalized-descent} follow from
\(\eta=\frac{\eta_0}{\sqrt T}\), \(C_1=\mathcal O(\eta)\), and
\(C_3,C_4=\mathcal O(\eta^2)\). Thus, \eqref{eq:app-normalized-descent}
gives
\[
    \frac1T
    \sum_{k=0}^{T-1}
    \mathbb E
    \left[
    \left\|
    \nabla\widetilde{\Phi}(u_k)
    \right\|^2
    \right]
    =
    \mathcal O\!\left(\frac1{\sqrt T}+
    \mathcal B^{(1)}
    +
    \mathcal B^{(2)}
    \right). \qedhere
\]
\end{proof}

\end{document}